\documentclass{article}

\usepackage[all]{xy}
\usepackage[dvips]{graphicx}
\usepackage{latexsym,amsthm,amssymb,amscd}
\usepackage{enumerate}
\usepackage{multicol}
\usepackage{amscd}
\usepackage{amsmath}
\usepackage{floatflt}

\numberwithin{equation}{section}
\numberwithin{figure}{subsection}
\begin{document}

\newcommand{\langerPfeil}{\begin{CD}@>>{hjkh}>\end{CD}}

\newcommand{\un}{\underline}
\newcommand{\ov}{\overline}

\newtheorem{theorem}{Theorem}[section]
\newtheorem*{theo}{Theorem}
\newtheorem*{conje}{Conjecture}
\newtheorem{lemma}[theorem]{Lemma}
\newtheorem{prop}[theorem]{Proposition}
\newtheorem{cor}[theorem]{Corollary}
\newtheorem{ex}[theorem]{Example}
\newtheorem{conjecture}[theorem]{Conjecture}
\newtheorem{definition}[theorem]{Definition}
\newcommand{\Bw}{{\it Proof:\quad}}
\newcommand{\Bwe}{\hfill $\diamondsuit$}

\newcommand{\cA}{{\mathcal A}}
\newcommand{\cB}{{\mathcal B}}
\newcommand{\cC}{{\mathcal C}}
\newcommand{\cD}{{\mathcal D}}
\newcommand{\cE}{{\mathcal E}}
\newcommand{\cF}{{\mathcal F}}
\newcommand{\cG}{{\mathcal G}}
\newcommand{\cH}{{\mathcal H}}
\newcommand{\cI}{{\mathcal I}}
\newcommand{\cK}{{\mathcal K}}
\newcommand{\cL}{{\mathcal L}}
\newcommand{\cM}{{\mathcal M}}
\newcommand{\cN}{{\mathcal N}}
\newcommand{\cO}{{\mathcal O}}
\newcommand{\cP}{{\mathcal P}}
\newcommand{\cS}{{\mathcal S}}
\newcommand{\cT}{{\mathcal T}}
\newcommand{\cU}{{\mathcal U}}
\newcommand{\cV}{{\mathcal V}}
\newcommand{\cW}{{\mathcal W}}
\newcommand{\cZ}{{\mathcal Z}}

\newcommand{\KL}{{Kazhdan-Lusztig}}
\newcommand{\Pol}{Kazhdan-Lusztig-polynomials}
\newcommand{\CF} {completion functor}

\newtheorem{Lemma}[theorem]{Lemma}
\newtheorem{korollar}[theorem]{Corollary}
\newtheorem{Prop}[theorem]{Proposition}
\newtheorem{de}[theorem]{Definition}
\newtheorem{Bed}[theorem]{Bedingung}
\newtheorem{Bez}[theorem]{Bezeichnungen}
\newtheorem{Bem}[theorem]{Bemerkung}
\newtheorem{remark}[theorem]{Remark}
\newtheorem{rk}[theorem]{Remark}
\newtheorem{Con}[theorem]{Convention}
\newtheorem{Verm}[theorem]{Vermutung}
\newtheorem{Bsp}[theorem]{Beispiel}
\newtheorem{Beob}[theorem]{Beobachtung}

\newcommand{\fC}{\mathbf{C}}
\newcommand{\fD}{\mathbf{K}}
\newcommand{\fP}{\mathbf{P}}
\newcommand{\fH}{\op{H}}
\newcommand{\mg}{\mathfrak{g}}
\newcommand{\mh}{\mathfrak{h}}
\newcommand{\mb}{\mathfrak{b}}
\newcommand{\mm}{\mathfrak{m}}
\newcommand{\mn}{\mathfrak{n}}
\newcommand{\ma}{\mathfrak{a}}
\newcommand{\mN}{\mathbb{N}}
\newcommand{\mK}{\mathbb{K}}
\newcommand{\mR}{\mathbb{R}}
\newcommand{\mV}{\mathbb{V}}
\newcommand{\mW}{\mathbb{W}}
\newcommand{\mC}{\mathbb{C}}
\newcommand{\mP}{\mathbb{P}}
\newcommand{\mZ}{\mathbb{Z}}
\newcommand{\SL}{\mathfrak{sl}}

\newcommand{\GO}{[\cO_0^\mZ]}

\newcommand{\HC}{Harish-Chandra}
\newcommand{\HI}{\cH^I}
\newcommand{\IH}{{^I}\!\cH}
\newcommand{\IHI}{{^I}\!\cH^I}
\newcommand{\IHJ}{{^I}\!\cH^I}
\newcommand{\JHI}{{^I}\!\cH^I}
\newcommand{\LAMU}{{_\la}\!\cH_\mu}
\newcommand{\LA}{{_\la}\!\cH_{\la'}}
\newcommand{\MU}{{_\mu}\!\cH_{\mu'}}
\newcommand{\Hoom}{{^{}_0\cH_0^{n}}}
\newcommand{\nHC}{{^n_0\cH_0^{}}}

\newcommand{\la}{\lambda}

\newcommand{\bs}{{\bf s}}
\newcommand{\be}{{\bf e}}

\newcommand{\oneone}{\stackrel{1:1}\leftrightarrow}
\newcommand{\ff}{\stackrel{f}\rightarrow}

\newcommand{\HOM}{\operatorname{Hom}}
\newcommand{\hOM}{\operatorname{hom}}
\newcommand{\gHOM}{\operatorname{gHom}}
\newcommand{\GHOM}{\operatorname{Hom}_\mg}
\newcommand{\EXT}{\operatorname{Ext}}
\newcommand{\eXT}{\operatorname{ext}}
\newcommand{\gEXT}{\operatorname{gEXT}}
\newcommand{\END}{\operatorname{End}}
\newcommand{\GEND}{\operatorname{End}_\mg}
\newcommand{\RES}{\operatorname{res}}
\newcommand{\p}{\mathfrak{p}}
\newcommand{\MOD}{\operatorname{mod}}
\newcommand{\MOF}{\operatorname{mof}}
\newcommand{\gMOD}{\operatorname{gmod}}
\newcommand{\gMOF}{\operatorname{gmof}}
\newcommand{\supp}{\operatorname{supp}}
\newcommand{\verg}{\operatorname{v}}
\newcommand{\op}{\operatorname}
\newcommand{\PDIM}{\op{pdim}}
\newcommand{\KRULL}{\op{Krull-Dim}}
\newcommand{\DEPTH}{\op{depth}}
\newcommand{\RANK}{\op{rank}}

\newcommand{\LLALA}{{_\la}L_\la}
\newcommand{\LLAMU}{{_\la}L_\mu}
\newcommand{\LLANU}{{_\la}L_\nu}
\newcommand{\LMUNU}{{_\mu}L_\nu}
\newcommand{\LNUMU}{{_\nu}L_\mu}
\newcommand{\LMUMU}{{_\mu}L_\mu}
\newcommand{\PLAMU}{{_\la}P_\mu}
\newcommand{\PMUMU}{{_\mu}P_\mu}
\newcommand{\PLANU}{{_\la}P_\nu}
\newcommand{\PMULA}{{_\mu}P_\la}
\newcommand{\PMUNU}{{_\mu}P_\nu}
\newcommand{\PNUMU}{{_\nu}P_{\mu}}
\newcommand{\PLAMUm}{\:^{}_\la P_\mu^n}
\newcommand{\PLAMUminus}{\:^{}_\la P_\mu^{m}}
\newcommand{\PLARHOm}{\:^{}_\laP_{-\rho}^n}
\newcommand{\PoRHOm}{{}^{}{_0}P_{-\rho}^n}
\newcommand{\Poom}{\:^{}_0P_{0}^n}
\newcommand{\PRHOom}{\:^{}_{-\rho}P_{0}^n}
\newcommand{\PRHORHOm}{\:^{}_{-\rho}P_{-\rho}^n}
\newcommand{\PLAMUml}{{^n_\la}P_\mu^{}}

\newcommand{\PLANUm}{\:^{}_\la P_\nu^n}
\newcommand{\PMULAm}{\:^{}_\mu P_\la^n}
\newcommand{\PMUNUm}{\:^{}_\mu P_\nu^n}
\newcommand{\PMUMUm}{\:^{}_\mu P_{\mu}^n}
\newcommand{\PNULAm}{\:^{}_{\nu}P_{\la}^n}

\newcommand{\PLAMUSm}{\:^{}{_{\la'}}P_{\mu'}^n}
\newcommand{\PLAMUSNm}{\:^{}_{\la'}P_{\mu}^n}

\newcommand{\PLAMUNSml}{\:{^n_{\la}}P_{\mu'}^{}}

\newcommand{\PLAMUnONE}{\:{_\la}^{n-1}P_\mu^{n-1}}
\newcommand{\PLANUnONE}{\:{_\la}^{n-1}P_\nu^{n-1}}
\newcommand{\PMUNUnONE}{\:{_\mu}^{n-1}P_\nu^{n-1}}
\newcommand{\PNUMUnONE}{\:{_\nu}^{n-1}P_\mu^{n-1}}

\newcommand{\PLAMUONE}{\:^{}_\la P_\mu^1}
\newcommand{\PLANUONE}{\:^{}_\la P_\nu^1}
\newcommand{\PMUNUONE}{\:^{}_\mu P_\nu^1}
\newcommand{\PNUMUONE}{\:^{}_\nu P_\mu^1}

\newcommand{\HLAMUn}{\:^{}_\la\cH_\mu^n}
\newcommand{\HLALAm}{\:^{}_\la\cH_\la^n}
\newcommand{\HLAMUeins}{\:^{}_\la\cH_\mu^1}
\newcommand{\HLAMUnl}{\:^{n}_\la\cH_\mu}
\newcommand{\HLAMUm}{\:^{}_\la\cH_\mu^n}
\newcommand{\HLANUn}{\:^{}_\la\cH_\nu^n}
\newcommand{\HMULAn}{\:^{}_\mu\cH_\la^n}
\newcommand{\HNUMUn}{\:^{}_\nu\cH_\mu^n}
\newcommand{\HMUNUn}{\:^{}_\mu\cH_\nu^n}
\newcommand{\HLAMU}{\:_\la\cH_\mu}
\newcommand{\HLANU}{\:_\la\cH_\nu}
\newcommand{\HMULA}{\:_\mu\cH_\la}
\newcommand{\HMUNU}{\:_\mu\cH_\nu}
\newcommand{\HNULA}{\:_\nu\cH_\la}
\newcommand{\HNUMU}{\:_\nu\cH_\mu}
\newcommand{\HMUNUmn}{\:^{}_\mu^m\cH_\nu^n}
\newcommand{\HNUMUmn}{\:^{}_\nu^m\cH_\mu^n}
\newcommand{\HoRHOm}{\:^{}_0\cH_{-\rho}^n}
\newcommand{\HoRHOeins}{\:^{}_0\cH_{-\rho}^1}
\newcommand{\HRHOom}{\:^{}_{-\rho}\cH_{0}^n}
\newcommand{\HRHOoeins}{\:^{}_{-\rho}\cH_{0}^1}
\newcommand{\HRHORHOm}{\:^{}_{-\rho}\cH_{-\rho}^n}

\newcommand{\TOR}{\operatorname{Tor}}
\newcommand{\HO}{\operatorname{H}}
\newcommand{\KER}{\operatorname{ker}}

\newcommand{\KERL}{{\kappa}}
\newcommand{\IDEAL}{\triangleleft}
\newcommand{\IM}{\operatorname{im}}
\newcommand{\COKER}{\operatorname{coker}}
\newcommand{\RAD}{\operatorname{rad}}
\newcommand{\ANN}{\operatorname{Ann}}
\newcommand{\RANN}{\operatorname{RAnn}}
\newcommand{\LANN}{\operatorname{LAnn}}
\newcommand{\DIM}{\operatorname{dim}}
\newcommand{\GK}{\operatorname{GKdim}\:}
\newcommand{\DHH}{\DIM\HOM_\cH}
\newcommand{\HH}{\HOM_\cH}
\newcommand{\DHS}{\DIM\HOM_{S\otimes S}}
\newcommand{\HS}{\HOM_{S\otimes S}}
\newcommand{\DHSW}{\DIM\HOM_{S^\mu\otimes S}}
\newcommand{\SOC}{\operatorname{soc}}
\newcommand{\DEG}{\operatorname{deg}}
\newcommand{\SUPP}{\operatorname{supp}}
\newcommand{\SPEC}{\operatorname{spec}}
\newcommand{\TR  }{\operatorname{Tr}}
\newcommand{\CX}{\operatorname{cx}}

\newcommand{\oast}{\circledast}

\newcommand{\ad}{^{\mbox\sl ad}}
\newcommand{\adf}{^{\mbox\sl adf}}
\newcommand{\opp}{^{\mbox\sl opp}}

\newcommand{\ON}{\theta_{(\la,\nu)}^{(\mu,\nu)}}
\newcommand{\OUT}{\theta_{(\mu,\nu)}^{(\la,\nu)}}
\newcommand{\RON}{\theta_{(\nu,\la)}^{(\nu,\mu)}}
\newcommand{\ROUT}{\theta_{(\nu,\mu)}^{(\nu,\la)}}
\newcommand{\On}{\theta_{on}}
\newcommand{\Out}{\theta_{out}}

\newcommand{\ZW}{\mathbb{Z}[\cW]}

\newcommand{\TM}{T}
\newcommand{\MM}{M_T(\lambda)}
\newcommand{\Pdef}{P_T(\lambda)}

\newcommand{\Ver}{M^i(\la)}
\newcommand{\Verma}{\cU\otimes_{\cU(\mb)}S(\mh)/(\ker \la)^i}
\newcommand{\NAB}{\nabla}
\renewcommand{\triangle}{\Delta}
\newcommand{\NO}{\nabla(0)}
\newcommand{\NxO}{\nabla(x\cdot\la)}
\newcommand{\VO}{\triangle(0)}
\newcommand{\VxO}{\triangle(x\cdot 0)}
\newcommand{\NABLA}{\nabla(\la)}
\newcommand{\Px}{P(x\cdot\la)}
\newcommand{\PxO}{P(x\cdot 0)}
\newcommand{\Lx}{L(x\cdot\la)}
\newcommand{\Ly}{L(y\cdot\la)}
\newcommand{\LxO}{L(x\cdot 0)}
\newcommand{\Mx}{M(x\cdot\la)}
\newcommand{\PANTI}{P(w_o\cdot\la)}
\newcommand{\MANTI}{M(w_o\cdot\la)}
\newcommand{\IANTI}{I(w_o\cdot\la)}

\newcommand{\T}{\operatorname{T}_\la^i}

\newcommand{\dann}{\Longrightarrow}
\newcommand{\dan}{\Rightarrow}
\newcommand{\inj}{\hookrightarrow}
\newcommand{\linj}{\longrightarrow}
\newcommand{\surj}{\mbox{$\rightarrow\!\!\!\!\!\rightarrow$}}
\newcommand{\iso}{\tilde{\rightarrow}}
\newcommand{\liso}{\tilde{\longrightarrow}}

\newcommand{\dota}{s_\alpha\cdot 0}
\newcommand{\dotb}{s_\beta\cdot 0}
\newcommand{\dotab}{s_\alpha s_\beta\cdot 0}
\newcommand{\dotba}{s_\beta s_\alpha\cdot 0}
\newcommand{\dotaba}{s_\alpha s_\beta s_\alpha\cdot 0}

\newcommand{\RZ}{$x=s_1\cdot\,\ldots\,\cdot s_r$}

\newcommand{\SIND}{\operatorname{S-Ind}}
\newcommand{\PO}{\operatorname{Pos}}
\newcommand{\CH}{\operatorname{ch}}
\newcommand{\CAN}{\operatorname{can}}
\newcommand{\ID}{\operatorname{id}}
\newcommand{\SYZ}{\operatorname{Syz}}

\newcommand{\lra}{\longrightarrow}
\newcommand{\pO}{\cO_0^\mathfrak{p}}
\newcommand{\OTL}{\cO_{T,\Lambda}}
\newcommand{\MTL}{\cM_{T,\Lambda}}
\newcommand{\Pinf}{P^\infty_T(\la)}
\newcommand{\PRinf}{P^\infty_{S(0)}(\la)}

\author{Catharina Stroppel}
\date{}
\title{TQFT with corners and tilting functors in the Kac-Moody case}
\maketitle

\begin{abstract}
We study projective functors (i.e. direct summands of compositions of
translations through walls) for parabolic versions of $\cO$ as well as for
integral regular blocks outside the critical hyperplanes in the symmetrizable Kac-Moody
case. It turns out that in both situations the functors are completely
determined by their restriction to the additive category generated by (the limit of) a `full projective tilting' object. We describe how
projective functors in the parabolic setup give rise to an invariant of tangle
cobordisms and formulate a conjectural direct connection to Khovanov homology. Our main result,
however, is the classification theorem for indecomposable projective functors
in the Kac-Moody case verifying a conjecture of F. Malikov and I. Frenkel.
\end{abstract}

\tableofcontents

\section*{Introduction}
The original motivation of this paper was the spectacular ``categorification
program'' described in \cite{BFK}, where the authors propose a way to get
tangle invariants via the representation theory of the Lie algebras
$\mathfrak{sl}_n$. A first step in this programme was carried out in
\cite{StDuke} where the main conjectures of \cite{BFK} were proved, thus providing a
functorial tangle invariant. The functors there are built up from the so-called projective functors by some cone
construction; these give tilting complexes and derived equivalences which
provide a functorial action of the braid groups.
The beautiful and amazing fact is that these functors lead to an
``enriched'' Jones invariant: the combinatorics of these
functors (on the level of the Grothendieck group) is given by the Jones polynomial for tangles (see \cite[Remark 7.3]{StDuke}). In the present paper we will mainly address the following three topics:

\subsubsection*{Invariants of tangles and cobordisms.}
We follow the philosophy of \cite{Rouquier} for example and consider not only
  the functors, which provide an invariant of tangles, but also natural
  transformations between them. This gives us
  (Theorem~\ref{2cat}, Theorem~\ref{2catorient}) a ``functorial realization''
  of the tangle $2$-category: The functorial invariant from \cite{StDuke}
  associates to each tangle diagram a functor leading to a tangle
  invariant (up to shifts). In the present paper we will show that the indeterminacy
  up to shifts disappears if we work with oriented tangles instead and get
  therefore a functorial invariant of oriented tangles. Moreover, the tangle
  invariants will be enriched by assigning to each cobordism between two
  oriented
  tangles a natural transformation between the corresponding functors. We will
  prove that, up to scalars, this defines an invariant of cobordisms. Hence,
  we construct a $3$-dimensional TQFT for ma\-ni\-folds with corner. It would be
  interesting to know whether the indeterminacy up to scalars can be explained
  in terms of an additional geometric or topological structure on the
  cobordisms.

\subsubsection*{Projective functors in the parabolic case and Khovanov homology.}
We show that projective functors (and their morphisms) for the
  parabolic versions of the Bernstein-Gelfand-Gelfand category $\cO$ are determined by evaluating at what we call a
  full projective tilting object (Theorem~\ref{Tcomb}). This generalises the
  classification result of \cite{BG} for the projective functors of the
  category $\cO$ to the parabolic versions of $\cO$, as well as the result of
  W. Soergel~\cite{Sperv} which says that the combinatorics of the category $\cO$
  is given by the endomorphism ring of a full projective tilting module. We
  conjecture that this result provides a direct link between the functorial
  invariants from \cite{StDuke} and the Khovanov homology introduced in
  \cite{Khovanov} and \cite{Khotangles} (Conjecture~\ref{Khovlink} and Example~\ref{itfits}). In particular, the combinatorics of Khovanov's homology should
  give a combinatorial description of the category of projective functors for
  certain parabolic categories $\cO$ corresponding to $\mathfrak{sl}_n$. On
  the other hand the approach towards invariants of tangles and cobordisms we
  will describe is much ``richer'' than the one of \cite{Khotangles}, since
  the underlying categories can be used to categorify (at least the main
  structures of) the tensor category of finite dimensional modules over
  the quantum group $\cU_q(\mathfrak{sl}_2)$.  It provides for example a
  categorical interpretation of the (dual) canonical bases and the Schur-Weyl
  duality. The first steps in this direction (using a Koszul dual setup) can
  be found in \cite{FKS}, based on \cite{BFK}.

\subsubsection*{The classification of projective functors in the Kac-Moody case:}
The classification theorem of projective functors in \cite{BG} and the
  Kazhdan-Lusztig combinatorics together can be formulated (in ``modern''
  language) as follows: {\it The projective endofunctors of the principal block of
  the category $\cO$ associated with a semisimple complex Lie algebra $\mg$
  form a ringoid which is isomorphic to the ringoid $\mZ[W]$
  of the associated Coxeter group.}
We will generalise this result to the case where
  $\mg$ is a symmetrizable Kac-Moody algebra. The first problem which occurs
  there is how to define projective functors.  In \cite{FM}, \cite{FMengl}
  such functors (so-called tilting functors) were defined for blocks of the
  category $\cO$ outside the critical hyperplanes associated to a
  symmetrizable Kac-Moody algebra. The definition uses the highly non-trivial
  Kazhdan-Lusztig tensoring. In \cite{FM} and \cite{FMengl}, I. Frenkel and
  F. Malikov also formulated a classification theorem in analogy to
  \cite{BG}. Unfortunately, there is a gap in the proof which, so far, cannot be
  fixed. Therefore, instead of following this path and using the
  Kazhdan-Lusztig tensoring, we define projective functors in terms of translation functors as they appear in \cite{Neidhardt1}, \cite{Neidhardt2} and
  \cite{CT}, and then prove the classification theorem proposed by I. Frenkel
  and F. Malikov.\\

We begin the paper by recalling the basics about parabolic category $\cO$. The
first result is Theorem~\ref{Tcomb} which gives a description of
homomorphism spaces between projective functors on the parabolic category
$\cO$. This result is probably well-known to specialists, but we were unable
to find it in the literature. The theorem states that the morphisms
between projective functors $F$ and $G$ are the bimodule
homomorphisms between $F(T)$ and $G(T)$, where $T$ is a full projective
tilting module. In the special case of the principal block $\cO_0$ of $\cO$ this comes
down to the well-known statement that morphisms between projective functors
are simply the bimodule morphisms between certain special bimodules defined by
W. Soergel. As an illustration we reprove the classification of projective
functors from~\cite{BG} using the Theorem~\ref{Tcomb} and Soergel's bimodules. Although the arguments are all taken from \cite{SHC}, there is a
slight (but for our purposes important) difference between our approach and
the ones in \cite{SHC} and \cite{BG}: they work either with
Harish-Chandra bimodules or with a deformed version of category $\cO_0$
with deformation ring being the (localised) universal enveloping algebra
$S(\mh)$; we work with a specialisation with respect to the centre $C$ of the category
$\cO_0$. For our purposes this is much stronger, since it can be generalised
to the case of Kac-Moody algebras where it is not clear how to define
Harish-Chandra bimodules. It fits also better with the
combinatorics of $\cO_0$ described in \cite{Sperv} and generalised in
\cite{Fiebigcomb} to blocks (outside the critical hyperplanes) of category $\cO$
for any symmetrizable Kac- Moody algebra.

This will be helpful in the last section when we define
projective functors for blocks (outside the critical hyperplanes) of symmetrizable Kac--Moody algebras based
on the translation functors defined in \cite{Neidhardt1}, \cite{Neidhardt2}
and \cite{CT}. By definition, a projective functor is nothing else than a
direct summand of direct sums of compositions of translations through
walls. This definition works over any deformation ring, in particular we will
deform with the centre $Z$ of the category. Any projective functor $F$ will
have a deformed version $F_Z$.

For the principal block $\cO_0$ of the classical category $\cO$ of some semisimple Lie
algebra $\mg$, the centre $Z$ is the
endomorphism ring of a full projective-tilting module $T$, namely the
indecomposable projective module with anti-dominant highest weight (\cite{Sperv}). On the
other hand, $Z$ is also the endomorphism ring of any $Z$-deformed Verma module
$\Delta_Z(\la)$ in the principal block (\cite{SHC}). As a special case of
Theorem~\ref{Tcomb} we get that any projective endofunctor $F$ is completely
determined by the $\mg\otimes_\mC Z$-module $F(T)$ which is isomorphic (as a
bimodule) to $F_Z\Delta_Z(0)$. The latter is a
projective object in the $Z$-deformed category $\cO$. If we specialise to the
usual category $\cO$ we get a projective object in the principal block of $\cO$. This
procedure gives a natural bijection between (indecomposable) projective
endofunctors of $\cO_0$ and (indecomposable) projective
objects in $\cO_0$ and a description of homomorphisms between such projective
functors (see Section~1).

For Kac-Moody algebras the principal idea will be the same, but there are
several obstacles to pass. Given a block $\cO_\Lambda$ (outside the critical
hyperplanes) it could happen that we have either (in the so-called positive level) enough projective modules, but no tilting
modules and no anti-dominant weight, or we have tilting modules and an
anti-dominant weight, but not necessarily projective modules (in the negative
level). In particular, there need not to be a full projective-tilting module. We first concentrate on blocks, where there is an
anti-dominant weight and use the idea of \cite{RW} (in the special situation
of \cite{CT}) that one should use a limit of projective modules from truncated
categories to get an analogue for Soergel's antidominant projective module and
then construct a fake full projective tilting module. We show in Theorem~\ref{centerKac} that its
endomorphism ring coincides with the centre $Z$ of the category $\cO_\Lambda$
(verifying a conjecture of P. Fiebig). Applying the tilting equivalence $\tau$ from \cite{SoKac}, which
is based on independent work of S. Arkhipov and A. Voronov, one obtains some
block $\cO_{\tau(\Lambda)}$ which contains a dominant weight $\tau(\la)$. There, we have the deformed Verma
module $\Delta_Z(\tau(\la))$ with endomorphism ring $Z$ and have also translation functors (through walls) as
defined in \cite{CT}. In particular, our previous definition of projective
endofunctors makes sense. Via the
equivalence $\tau$ any such projective functor $F$ induces some endofunctor,
say $\tau(F)$, of the block we started with. We call these functors also projective. Via $\tau$, the
projective modules in blocks of positive level become tilting modules in
blocks of negative level, hence these projective functors $\tau(F)$ map
tilting modules to tilting modules. This is the reason why what we call ``projective functors'' are
called ``tilting functors'' in \cite{FMengl}. In analogy to the classical
situation we will prove the following result which is a natural generalisation
of \cite{BG} and coincides with the conjectural classification theorem in
\cite{FM} (note the typo in the formulation there).
\begin{theo}
Let $\mg$ be a symmetrizable Kac-Moody algebra over $\mC$. Let
$\cO_{\mC,\Lambda}$ be a regular block, outside the critical hyperplanes, of the corresponding category $\cO$
such that $\Lambda$ contains an antidominant weight $\lambda$.
\begin{enumerate}[(a)]
\item Let $Z$ denote the centre of $\cO_{\mC,\Lambda}$. Let $F_\mC, G_\mC:\cO_{\mC,\Lambda}\rightarrow\cO_{\mC,\Lambda}$ be compositions of translations through
walls with corresponding deformed functors
$F_{Z}, G_Z:\cO_{Z,\Lambda}\rightarrow\cO_{Z,\Lambda}$. Let $\Delta_Z(\la)$ be the
$Z$-deformed Verma module with highest weight $\la$. Then there is an isomorphism of vector spaces (or even of rings if $F_\mC\cong G_\mC$)
\begin{eqnarray*}
\HOM(F_\mC,G_\mC)\cong\HOM_{\mg\otimes Z}(F_{Z}\Delta_{Z}(\la),G_{Z}\Delta_{Z}(\la)).
\end{eqnarray*}
\item
  There are natural bijections of isomorphism classes
\begin{eqnarray*}
  \begin{array}[t]{cc}
    \big\{\text{indecomposable projective endofunctors of } \cO_{\mC,\Lambda})\}&F\\
\updownarrow \text{\small 1:1}&\downarrow\\
\big\{\text{indecomposable tilting objects of }\cO_{\mC,\Lambda}\big\}&F(\Delta_\mC(\la))\\
\updownarrow \text{\small 1:1}&\downarrow\\
\big\{\text{indecomposable projective objects of }\cO_{\mC,\tau(\Lambda)}\big\}&\tau(F(\Delta_\mC(\la)))\\
\updownarrow \text{\small 1:1}&\downarrow\\
 \big\{\text{indecomposable projective endofunctors of } \cO_{\mC,\tau(\Lambda)}\}&\tau(F)
  \end{array}
\end{eqnarray*}
\end{enumerate}
\end{theo}

The theorem implies in the classical case of a finite dimensional Lie algebra
the results of \cite[Section 4]{Erikproj}. With the results from
Kazhdan-Lusztig theory, it follows that {\it the projective functors from
  $\cO_\Lambda$ to itself categorify the group algebra of the Coxeter group
  corresponding to the block $\cO_\Lambda$.}\\

Certainly, one would like to have an explicit description of the centre $Z$. A
first step to solve this problem was done in \cite{CT}, where the centre of a certain deformation $\cO_{S_{(0)},\Lambda}$ of
  $\cO_{\mC,\Lambda}$ was determined. However, it remained an essential open
  question if specialisation gives rise to the specialised centre (as was
  implicitly claimed in \cite{CT}). The methods of \cite{CT} are not
  sufficient to confirm this claim. In general, the deformed situation is quite
  different from and much "easier" than the specialised one (see
  e.g. \cite[Theorem 3.12]{Fiebigcomb} in comparison to Soergel's Structure
  Theorem \cite{Sperv}). However, we prove in Theorem~\ref{centerKac} that
  the centre behaves well under specialisation, and together with the main
  result of \cite{CT} it turns out that the endomorphism ring of the
  fake full projective module mentioned above is a completion of the cohomology ring (see e.g. \cite{KacKaz}) of the
full flag manifold corresponding to the Langlands dual group. In particular, it has a natural grading. Therefore, the next step
would be to consider a graded version of the category $\cO$ for Kac--Moody
algebras in analogy to \cite{BGS} and define graded versions of projective
functors using the approach of \cite{StGrad}. In this
way one should get a categorification of the generic Hecke algebra instead of
just a categorification of the group algebra.\\

Let us come back to our invariants of tangles and cobordisms mentioned
earlier. A classical knot is an embedded circle in $\mR^3$. Knot diagrams are
generic projections of the image of such an embedding onto the plane. Two knot diagrams represent (up to ambient isotopy) the same
knot if and only if they are related via a sequence of Reidemeister
moves. These moves can be considered as surfaces in $\mR^3\times [0,1]$ properly mapped into $\mR^2\times [0,1]$, where the
boundaries  $\mR^2\times \{0\}$ and  $\mR^2\times \{1\}$ are strings before and
after the move. On the other hand, any such knot diagram can be decomposed
into generating elementary parts by introducing a so-called height function. We have (see e.g. \cite[Theorem 2.2.1]{CRS}) the following
\begin{enumerate}[Fact(R):]
\item two knot diagrams with a height function represent (up to ambient isotopy) the same
knot if and only if they are related via a sequence of Reidemeister,
$T$-, $H$- and $N$-moves. (see Section~\ref{inv}, Figure~\ref{fig:2morphisms}).
\end{enumerate}

Instead of working with knots we will consider tangles. That means we fix a
finite number (say $m$ and $n$) of points in  $\mR^2\times \{0\}$ and
$\mR^2\times \{1\}$ respectively. An $(m,n)$-tangle is a set of disjoint smooth
curves in  $\mR^2\times [0,1]$, intersecting the boundary in $n+m$ points,
which are exactly the $n+m$ endpoints of all the non-closed smooth curves. Any
generic projection to the plane gives a tangle diagram. As for knot diagrams,
Fact (R) holds and
all the tangles can be written as a composition of elementary tangles which
are depicted in Figure~\ref{fig:elementary}. Just as for knots, the
moves have interpretations in terms of generic (knotted) surfaces (or cobordisms) properly embedded into
$\mR^3\times [0,1]$, where the boundaries are the tangles before and after the
move. Such a surface can be described via a family $D_t$, $t\in [0,1]$ of tangle
diagrams by first
projecting the surface to $\mR^2\times [0,1]$, and then letting $D_t$ be the
diagram representing the intersection with the planes
$\mR^2\times\{t\}$. There are only finitely many critical points for $t$,
which  means between these points the diagrams undergo just planar
isotopies. One can therefore describe the surfaces by a finite sequence of diagrams which
represent the intersections at critical points. These sequences are called
movies. In this way, any cobordism (or knotted surface in \cite{CRS}) between
tangles can be described via a movie. In Figure~\ref{fig:moves1114} one can
find $8$ movies displayed. D. Roseman~\cite{Roseman} defined Reidemeister
moves for surfaces and proved that two movies represent (up to ambient
isotopy) the same knotted surface if and only if they are related via a
sequence of Reidemeister moves.

In \cite{CRS} the authors gave a combinatorial description in form of a list of {\it elementary string interactions} and proved
that the diagrams associated to any knotted surface is a finite sequence of
elementary string interactions and that to each such finite sequence there exists
in fact a corresponding knotted surface. Some elementary string interactions
are depicted in Figure~\ref{fig:moves1114}.
The main theorem of \cite{CRS} gives
a list of {\it movie moves} such that {\it two sequences of elementary string interactions represent (up to ambient isotopy) the same
knotted surface if and only if they are related via a sequence of movie
moves}. Examples of such movie moves are depicted in Figure~\ref{fig:moves1114}.\\

In Theorem~\ref{2catorient} we will extend the
functorial invariant of tangles from \cite{StDuke} to an invariant of (oriented) tangles and
cobordisms. To any tangle diagram we associate a functor and to each elementary string
interaction we associate a natural transformation and show that these natural
transformations satisfy (up to scalars) the relations from \cite{CRS}. Hence
we get as a result a functorial invariant of oriented tangles and cobordisms
which can be formulated as follows (using the notation from Section~\ref{Section1} and Section~\ref{inv}):

\begin{theo}
  There is a $2$-functor $\Phi^{or}:{\mathcal Tan}^{or}\rightarrow{\mathcal
  Func}$ such that
\begin{enumerate}
\item if $t_1$ and $t_2$ are $1$-morphisms which differ by a sequence of
  Reidemeister, {T-}, {H-} or N-moves then there in isomorphism of functors
  $\Phi^{or}(t_1)\cong\Phi^{or}(t_2)$.
\item if $c_1$ and $c_2$ are sequences of generating $2$-morphisms which
  differ by a sequence of movie moves then $\Phi^{or}(c_1)=\Phi^{or}(c_2)$.
\end{enumerate}
\end{theo}

Invariants of tangles and cobordisms were already
obtained by M. Jacobsson and M. Khovanov (\cite{Jacobsson},
\cite{Khovanovcob}). These two papers contain an ``enrichment'' of the knot and tangle invariant introduced by M. Khovanov (\cite{Khovanov}, \cite{Khotangles}) who assigned to each
tangle or link diagram the homology of a combinatorially defined complex of
$\cH_n$-bimodules for some explicit given algebra $\cH_n$. (For a very nice
overview with simplified arguments we refer to \cite{BarNatan}).\\

In the first sight, there is no connection between the homology introduced by
M. Khovanov and the approach proposed in \cite{BFK} and \cite{StDuke}. However, it turns out that
the arguments which establish the extension of \cite{StDuke} to an invariant of
cobordisms just mimic the arguments given in \cite{Khovanovcob}, provided one
actually proves that the functors have similar properties to the functors
given by tensoring with the complexes of $\cH_n$-bimodules considered in
\cite{Khotangles}. In Section~\ref{inv} we will establish this and show
that the functors share all the important nice properties. In \cite{Braden} it is already mentioned that a certain
parabolic category $\cO$ for $\mathfrak{sl}_n$ is equivalent to a module
category over an algebra $A$ such that the algebra $\cH_n$ used in
\cite{Khotangles} is a subquotient of $A$. The algebra $A$ is described in terms
of quivers and relations in \cite{Braden}. Built on the results of T. Braden, we will
formulate a conjecture (Conjecture~\ref{Khovlink}) which bridges the approaches
of \cite{BFK} and \cite{Khovanov}. It says the following.
Based on \cite{BFK}, we associate in particular to each
$(2m,2n)$-tangle $t$ a functor
$\Phi^{or}(t):\bigoplus_{k=0}^{2n}{\bf
  D}^b(\gMOF\text{-}A_{2n}^k)\rightarrow\bigoplus_{k=0}^{2m}{\bf D}^b(\gMOF\text{-}A_{2m}^k)$,
where ${\bf D}^b(\gMOF\text{-}A_n^k)$ denotes the bounded derived category of the graded
version of $\cO^{\p_k}_0(\mathfrak{sl}_{2n})$ which is the principal block of the
 parabolic category $\cO$ corresponding to $\mathfrak{sl}_{2n}$ with
 parabolic Weyl group isomorphic to $S_k\times S_{2n-k}$. Recall that for each
 $\cO_0^{\p_k}(\mathfrak{sl}_{2n})$ we have a full projective tilting module
 $T_{2n}^k$. If we first restrict
 $\Phi^{or}(t)$ to a functor from ${\bf D}^b(\gMOF\text{-}A_{2n}^n)$ to
 ${\bf D}^b(\gMOF\text{-}A_{2m}^k)$, and then to perfect complexes of projective tilting
modules we finally assign to each $(2m,2n)$-tangle a functor
which can be realized as tensoring with a complex $\check{X}(\Phi^{or}(t))$ of
$(\END_\mg(T_{2m}),\END_\mg(T_{2n}))$-bimodules. We conjecture the following
direct connection to Khovanov homology:

\begin{conje}
  \begin{enumerate}
  \item For any natural number $m$, there is an isomorphism of algebras
  $p_m:\END_\mg(T_{2m})\cong\cH_m$, where $\cH_m$ denotes Khovanov's algebra.
  \item The homological tangle invariant $t\mapsto {\bf
  H}^\bullet(\check{X}_{\Phi^{or}(t)})$ is Khovanov's invariant.
  \end{enumerate}
\end{conje}

This  conjecture will be illustrated for $n=m=2$ in Example~\ref{itfits}. If this conjecture is true, we get Khovanov's tangle homology as a
special case of our much more general approach. Given such a direct connection would embed Khovanov's approach into the much
richer structure coming from \cite{BFK}, where most of the monoidal category
of (quantised) $\mathfrak{sl}_2$-modules can be seen (see for example
\cite{FKS}), and vice versa it would provide the beautiful, computable and
understandable combinatorics from Khovanov homology to describe the functorial tangle invariants.

\subsubsection*{Acknowledgements}
I would like to thank Henning Haahr Andersen, Peter Fiebig, Igor Frenkel, Volodymyr Mazorchuk and
Wolfgang Soergel for many helpful discussion related to this work. I
also thank Tom Braden for detailed explanations and information about his work
and Iain Gordon for comments on a previous version of the paper. I
finally would like to thank the organisers of the BIRS meeting
``Representations of Kac--Moody algebras and combinatorics'' and the
organisers of the ARTIN meeting 2004 who gave me
the opportunity to present parts of these results. I am grateful to EPSRC for
supporting me as a research associate.


\section{The finite dimensional parabolic situation}
\label{Section1}
In this section we study the principal block of parabolic versions of the
Bernstein-Gelfand-Gelfand category $\cO$. We define the notion of projective
functors generalising the definition in \cite{BG} and show that they are
determined already by their restriction to projective modules which are also
injective or tilting. The main result will be a theoretical description of morphism
spaces between projective functors. It turns out that they are naturally
graded. With the correct choice of a parabolic in type $A$, conjecturally the
endomorphism ring of projective-injective modules gives rise to the algebra
which defines Khovanov's homology and therefore provides a representation
theoretic interpretation of Khovanov's combinatorial approach.

\subsection{Notations and Preliminaries}
For any ring $R$ we denote by 
$\MOD\text{-}R$ ($\MOF\text{-}R$) and $R\text{-}\MOD$ ($R\text{-}\MOF$)
respectively the category of (finitely generated) right/left $R$-modules. Likewise, $R\text{-}\MOF\text{-}S=R\otimes
S^{\op{opp}}\text{-}\MOF=\MOF\text{-}R^{\op{opp}}\otimes S$ denotes the category of finitely
generated $(R,S)$-bimodules. If, additionally, $R$ and $S$ are graded rings
the symbols $R\text{-}\gMOF$, $\gMOF\text{-}R$, $R\text{-}\gMOF\text{-}S$  etc. denote the corresponding
categories of {\it graded} modules with degree preserving morphisms. In the following {\it graded} always means $\mZ$-graded. If
$M=\oplus_{i\in\mZ} M_i\in\gMOF\text{-}R$ is a graded module and $n\in\mN$, we denote by $M\langle
n\rangle$ the object in $\gMOF\text{-}R$ such that $M\langle n\rangle =M$ as right
$R$-module, but the grading shifted by $n$, i.e. $(M\langle
n\rangle)_i=M_{i-n}$. By an {\it algebra} we always mean a finite dimensional
unitary associative algebra over the complex numbers.

Let $\mg$ be a finite dimensional semisimple complex Lie algebra with its universal enveloping
algebra $\cU(\mg)$. We fix $\mb\supset\mh$, a Borel and a Cartan
subalgebra. Let $\cU(\mb)$ and $\cU(\mh)$ be the corresponding universal enveloping algebras. Let $\cO=\cO(\mg, \mb)$ denote the corresponding
highest weight category defined by Bernstein, Gelfand and Gelfand (\cite{BGG}). The objects
of this category are finitely generated $\cU(\mg)$-modules $M$ which are
locally $\cU(\mb)$-finite and have a weight space decomposition $M=\oplus
M_\la$, where $M_\la=\{m\in M\mid hm=\la(h)m\}$ for any $h\in\mh$. For
precise definitions and properties of $\cO$ we refer for example to \cite{BG}, \cite{Ja1}, \cite{Ja2}.

Let $W$ be the Weyl group corresponding to $\mg$ with longest element $w_0\in
W$ and let $\rho$ denote the
half-sum of positive roots. The category $\cO$ has a decomposition into
blocks, more precisely $\cO=\bigoplus\cO_{\la-\rho}$, where $\la$ runs through
the set of dominant weights and $\cO_{\la-\rho}$ denotes the block containing
the simple module with highest weight $\la-\rho$. In particular, $\cO_{-\rho}$
is semisimple and $\cO_0$ is the principal block containing the trivial
representation.

 We fix a parabolic subalgebra $\p$ containing $\mb$ and consider the parabolic category $\cO^\p$ defined as the full
subcategory of $\cO$ given by locally $\p$-finite objects. We refer to \cite{RC} and
\cite{Irself} for properties. We {\bf fix} $P^\p$, a minimal projective generator of the
main block $\cO_0^\p$ of $\cO^\p$ and denote
$A^\p=\END_\mg(P^\p)$ its endomorphism ring. The block decomposition of $\cO$
induces a decomposition  $\cO^\p=\bigoplus\cO^\p_{\la-\rho}$. Note that the
direct summands could become decomposable or even trivial. This is however not
 the case for the principal block $\cO^\p_0$ on which we will focus our attention.

\subsection{The algebras $A^\p=\END_\mg(P^\p)$}
Since $P^\p$ is a projective generator, the functor
$\HOM_\mg(P^\p,\bullet)$ defines an equivalence of categories
\begin{eqnarray}
\label{equiv}
  \epsilon^\p:\cO_0^\p\iso\MOF\text{-}A^\p,
\end{eqnarray}
where $\MOF\text{-}A^\p$ denotes the category of finitely generated right
$A^\p$-modules (see e.g. \cite[Section 2]{Bass}).

\begin{ex}
\label{ex}{\rm
 For $\mg=\mathfrak{sl}_n$ and $\p$ a maximal parabolic subalgebra with
  corresponding Weyl
  group isomorphic to $S_{1}\times S_{n-1}$, the algebra $A^\p$ is isomorphic
  to the path algebra of the quiver
  \begin{displaymath}
    \xymatrix{{\bullet}
\ar@/^/[r]^{a_1}&\bullet\ar@/^/[l]^{b_1}\ar@/^/[r]^{a_2}&\bullet\ar@/^/[l]^{b_2}\ar@{..}[r]&\bullet\ar@/^/[r]^{a_{n-1}}&\bullet\ar@/^/[l]^{b_{n-1}}
}
  \end{displaymath}
with vertices $1, 2, \ldots, n$ and the $2n$ arrows as indicated, with the relations $b_1a_1=0$,
$a_{i+1}a_i=0=b_{i}b_{i+1}$ and $a_ib_i=b_{i+1}a_{i+1}$ whenever the
expression makes sense. This is exactly the algebra appearing in \cite{KS}.
}
\end{ex}

In general, a (handy) explicit description of the algebra $A^\p$ is not
  known. However one can find partial information in the literature. Let us
  mention some: For type $A$ or $D$ and $\p$ a maximal parabolic, the paper
  \cite{Braden} of T. Braden gives a very nice and
  useful, but unfortunately not very handy, description of the algebras $A^\p$
  in terms of quivers with relations. For the non-parabolic case $\p=\mb$ some
  explicit examples can be found in \cite{Stquiv}, obtained by an algorithm
  based on \cite{Sperv} and recently generalised and substantially improved in \cite{Vybornov}. All the algebras $A^\p$ can be equipped with a
  grading which turns them into Koszul algebras (\cite{BGS}). In the
  Example~\ref{ex} the grading is obtained by putting every path in degree
  one. The dimension of the algebra can be determined using the combinatorics of parabolic Kazhdan-Lusztig polynomials (namely the $n_{x,y}$
   in the notation of \cite{SoKipp}; see e.g. \cite[Theorem 3.11.4]{BGS}). The
  representation type of the algebras in question is determined in
  \cite{BN}. The categories $\cO^\p$ are highest weight categories (in
the sense of \cite{CPS}), hence the algebras $A^\p$ are quasi-hereditary in the sense of \cite[Appendix]{Donkin}.

\subsection{$\cO_0^\p$ as highest weight category}
For $\la\in\mh^*$ let $W_\la=\{w\in W\mid w\cdot\la\}$ denote the stabiliser
of $\la$ under the dot-action $w\cdot\la=w(\la+\rho)-\rho$. We denote by $W_\p\subseteq W$ the parabolic subgroup corresponding to $\p$.
Let $W^\p$ denote the set of minimal coset representatives in
$W_\p\backslash W$. The category $\cO_0^\p$ is a highest weight category where the {\it standard objects} are the
(generalised) Verma modules $\Delta^\p(x\cdot0)$ with highest weights
$x\cdot0$ for $x\in W^\p$. Recall that $\Delta^\p (x\cdot0)$ is the maximal
quotient, contained in $\cO^\p$, of the Verma module $\Delta(x\cdot0)=\cU(\mg)\otimes_{\cU(\mb)}\mC_{x\cdot0}$. We denote by $L(x\cdot0)=L^\p(x\cdot0)$ the simple
head of $\triangle^\p(x\cdot0)$ and by $P^\p(x\cdot0)$ its projective
cover. (The latter is always considered as a direct summand of $P^\p$.) By
abuse of language, the standard, simple and projective modules in $\MOF\text{-}A^\p$ are
denoted by the same symbols. A module in $\cO_0^\p$ (or in $\MOF\text{-}A^\p$
respectively) {\it has a standard flag}, if it has
a filtration with subquotients isomorphic to generalised Verma modules (or
standard objects). Any projective module has a standard flag.
\subsection{The full projective tilting module}
Let $T=T^\p\in\MOF\text{-}A^\p$ be a {\it full projective tilting} module, i.e. the
direct sum over all modules constituting a system of representatives for the
isomorphism classes of indecomposable objects which are at the same time
projective, tilting and injective in $\MOF\text{-}A^\p$. (Note that two of these
three properties automatically force the third one, since there is a
contravariant duality which preserves simple objects.)
In \cite{Irself}, Irving studied projective modules in $\cO_0^\p$ and found
a very nice characterisation of projective modules which are also
injective. We summarise the main results of his paper in the following
proposition:
\begin{prop}(\cite{Irself})
\label{Irving}
  Let $P\in\cO_0^\p$ be projective. Then the following are equivalent
  \begin{enumerate}[(i)]
  \item $P$ is injective in $\cO_0^\p$,
  \item $P$ is tilting in $\cO_0^\p$,
  \item the head of $P$ contains only simple modules of maximal possible
  Gelfand-Kirillov dimension,
  \item any composition factor of the head of $P$ occurs as a submodule in
  some standard module $\Delta^\p(x\cdot0)$,
  \item $P$ is isomorphic to some direct sum of injective hulls of standard modules.
  \end{enumerate}
\end{prop}

We choose $T$ as a submodule of $P^\p$ and its endomorphism ring
$D^\p=\END_\mg(T)$ as a subalgebra of $A^\p$. Note that $D^\p$ is always a
finite dimensional Frobenius algebra. If we are in the special situation, where
$\p=\mb$ then $T$ is the unique indecomposable projective-injective module (equivalent to the
projective cover of the simple Verma module) with endomorphism ring
$C=\cU(\mh)/({\cU(\mh)}_+^W)$, the algebra of coinvariants (\cite{Sperv}). In
particular, $D^\mb$ is a commutative symmetric algebra.

\begin{ex}
\label{ex2}
  {\rm Let us consider the Example~\ref{ex} and denote the indecomposable
  projective module corresponding to vertex $i$ by $P(i)$.
    \begin{enumerate}[a.)]
    \item \label{exa}
 If $n=3$ then the vertices $1$, $2$, $3$ correspond to the simple objects $L(0)$,
    $L(s\cdot0)$ and $L(st\cdot0)$, where $s$, $t$ are simple reflections. The socle
    series of the indecomposable projective modules are of the form
 $P^\p(0)=L(0)|L(s\cdot0)$, $P^\p(s\cdot0)=L(s\cdot0)|L(0)\oplus L(st\cdot0)|L(s\cdot0)$
 and $P^\p(st\cdot0)=L(st\cdot0)|L(s\cdot0)|L(st\cdot0)$. Note that
 $P^\p(0)=\Delta^\p(0)$, $\Delta^\p(s\cdot0)=L(s\cdot0)|L(st\cdot0)$ and
 $\Delta^\p(st\cdot0)=L(st\cdot0)$. We have $T=P(2)\oplus P(3)$. For arbitrary
 $n\geq 2$, we have $T=\bigoplus_{i=2}^n P(i)$. If $n=2$ then $T=P(2)$ and its endomorphism algebra is the graded algebra
  $\mC[x]/(x^2)$, where $x$ is of degree $2$. This is exactly the algebra $\cA\langle 1\rangle$ in \cite{Khovanov}.
\end{enumerate}
}
\end{ex}

In general, the
module $T$ is quite difficult to describe and its endomorphism algebra is not
known. It is known, however, that for $\mg=\mathfrak{sl}_n$ this algebra is
symmetric, and depends only on the composition describing $\p$, but not on
the partition (see \cite{MSsym}). Conjecturally the centre of this algebra is
the cohomology ring of the associated Springer fibre (see \cite{KhSpringer}).

\subsection{The structure theorem for the category $\cO_0^\p$}
The philosophy behind our approach and one reason why we want to consider full
projective-tilting modules is that from the knowledge of the full-tilting module $T$ with
its endomorphism ring $\END_\mg(T)$ one could in principle recover the whole
category $\cO^\p$. This point of view fits perfectly well with the special
case of Soergel's description of $\cO_0^\mb$ in terms of modules over the
coinvariant algebra $\END_\mg(T^\mb)$.

An important property of the full tilting module is given by the following proposition whose proof relies on the fact that projective objects in
$\cO_0^\p$ can be built up from the projective Verma module using translation
functors and makes clear how one should understand Irving's result from
Proposition~\ref{Irving}. It naturally generalises Soergel's structure theorem
(\cite{Sperv}). To formulate it we need a little bit more notation. Let $\theta_s:\cO_0\lra\cO_0$ denote the {\it translation through}
the $s$-wall as defined for example in \cite[Section 3]{GJ}. The functor
$\theta_s$ is exact and its own biadjoint, hence a so-called Frobenius
functor. Note that $\theta_s$ maps the Verma module $\Delta(0)$ to the indecomposable projective module
$P(s\cdot0)$. By abuse of language we denote by the same symbol also its
restriction $\theta_s:\cO_0^\p\lra\cO_0^\p$ as well as the induced endofunctor
of $\MOF\text{-}A^\p$ (via \eqref{equiv}).

\begin{definition}{\rm
We call a functor $F:\cO_0^\p\lra\cO_0^\p$ (or
$F:\MOF\text{-}A^\p\lra\MOF\text{-}A^\p$) {\it
  projective} if it is a direct sum of direct summands of some composition of
translations through walls. In this case we also say $F$ {\it is a projective
functor on} $\cO_0^\p$ or $\MOF\text{-}A^\p$ respectively.}
\end{definition}

\begin{remark}{\rm
Note, that the classification theorem of projective functors from \cite{BG}
 (see Theorem~\ref{class} below) implies that our usage of the notation {\it
 projective functor} on $\cO^\mb_0$ is compatible
with the definition in \cite{BG}. In this
classification the functor
$\theta_s$ is then the (up to isomorphism unique) projective
functor mapping $\Delta(0)$ to the projective module
$P(s\cdot0)$. Since $P(s\cdot0)$ is indecomposable, so is the functor
 $\theta_s$ (for a
direct proof see Lemma~\ref{indec} below).}
\end{remark}

In the Example~\ref{ex2}~\eqref{exa}, we have an inclusion of $P(0)$ into
$P(s\cdot0)$ such that the cokernel can be embedded into $P(st\cdot0)$. In
general the following holds:
\begin{Prop}
  \label{copres}
  Let $P\in\gMOF\text{-}A^\p$ be projective. Then there exists a {\it
  projective-tilting-copresentation}, that is an exact sequence of the form
\begin{eqnarray*}
  0\rightarrow P\lra\bigoplus_I T\lra\bigoplus_J T,
\end{eqnarray*}
for some finite index sets $I$, $J$.
\end{Prop}
\begin{proof}
  By Proposition~\ref{Irving}, a module from $\pO$ is indecomposable
  projective-tilting
  if its socle (=head) is contained in the socle of a Verma
  module; and visa versa the projective hull of any simple composition factor
  occurring in the socle of a Verma module is injective. Therefore, any standard
  module, or even any module from $\gMOF\text{-}A^\p$ having a standard flag, embeds into a finite direct sum
  of copies of $T$. If $P=\Delta^\p(0)$
  is the projective standard module, then (by weight considerations) the cokernel of this
  embedding has again a standard filtration. Therefore, the statement of the lemma is
  true for $P=\Delta^p(0)$. Hence, there is a copresentation
\begin{eqnarray*}
  0\rightarrow \Delta^\p(0)\stackrel{i}\lra\bigoplus_I T\lra\bigoplus_J T,
\end{eqnarray*}
for some finite sets $I$ and $J$. Now, any projective module is of the form $F\Delta^\p(0)$ for some projective
functor $F$ (\cite[Proposition (v)]{Irself}). On the other hand, the socles of
the indecomposable projective-injective modules in $\cO_0^\p$ are exactly the
simple objects with maximal Gelfand-Kirillov dimension (this is Proposition~\ref{Irving}). Recall
the following general fact: Assume $M\in\cO_0^\p$ and
$\HOM_\mg(L,M)=0$ for any simple object in $\cO_0^\p$ not having maximal
Gelfand-Kirillov dimension and let $F:\cO_0^\p\rightarrow\cO_0^\p$ be a
projective functor. Then, by definition, $FM\in\cO_0^\p$, but we also claim
that $\HOM_\mg(L,FM)=0$ for any simple object in $\cO_0^\p$ not having maximal
Gelfand-Kirillov dimension. To show this let $G$ be the adjoint functor of
$F$. This is again a projective functor and we get
$\HOM_\mg(L,FM)=\HOM(GL,M)=0$, since $G$ does not increase the
Gelfand-Kirillov dimension (\cite[Lemma 8.8]{Ja2}), and therefore any
quotient of $GL$ has smaller Gelfand-Kirillov dimension than any arbitrary
non-zero submodule of $M$ (\cite[Lemma 8.6]{Ja2}). The claim follows.\\
In particular, the socle of the cokernel $K$ of $F(i)$ has only
composition factors of maximal Gelfand-Kirillov dimension. Applying again the
Proposition~\ref{Irving} we know that the injective hull of $K$ is also
projective, hence a direct summand of some $\bigoplus_{I'}T$. On the other
hand $F(T)$ is also projective and injective, hence a  direct summand of some
$\bigoplus_{I'}T$. The statement for general $P$ follows.
\end{proof}

The following generalises \cite[Struktursatz]{Sperv} (see \cite[Theorem 10.1]{Stquiv})
\begin{cor}
\label{faithful}
  The functor $\mV^\p=\HOM_{\MOF\text{-}A^\p}(T,\epsilon(\bullet)):\cO_0^\p\rightarrow\MOF\text{-}D^\p$ is fully
  faithful on projectives, i.e. it induces a natural isomorphism
  \begin{eqnarray*}
    \HOM_\mg(P_1,P_2)\cong\HOM_{\MOF\text{-}D^\p}(\mV^\p P_1,\mV^\p P_2)
  \end{eqnarray*}
for projective objects $P_1$, $P_2\in\cO_0^\p$.
\end{cor}
\begin{proof}
  Since any simple object occurring in the socle of a projective object in
  $\cO_0^\p$ is not annihilated by $\mV^\p$ (Proposition~\ref{Irving}), the natural map
  is injective.  It is obviously an isomorphism if $P_2=T$, since the dimension on both sides is just the number of simple composition factors of $P_1$ of
  maximal Gelfand dimension which equals the dimension of $\mV^\p P_1$. Then it
  is also an isomorphism if $P_2$ is a finite direct sum of copies
  of $T$. For the general case, we take a copresentation $\cC$ of $P_2$ coming from an
  injective-tilting copresentation for $\epsilon{P_2}$ via $\epsilon^{-1}$. Then
  $\mV^\p\cC$ is exact and stays exact when applying
  $\HOM_{\MOF\text{-}D^\p}(\mV^\p P_1,\bullet)$. On the other hand, applying
  $\HOM_\mg(P_1,\bullet)$ to $\cC$ is also exact, since $P_1$ is
  projective. The desired result follows then using the Five lemma.
\end{proof}

\subsection{Projective functors and projective tilting objects}
Any projective functor from $\MOF\text{-}A^\p$ to $\MOF\text{-}A^\p$ preserves the additive
category of projective tilting objects. In the following we want to show that
the functor is already determined by its restriction to this additive
category.

Let $F$ be a projective functor from $\MOF\text{-}A^\p$ to $\MOF\text{-}A^\p$ then $F(T)\in\MOF\text{-}A^\p$ by
definition. It has also a left $D^\p$-module structure given by $d.t=F(d)(t)$ for
any $d\in D^\p$, $t\in F(T)$ giving rise to a $(D^\p,A^\p)$-bimodule
structure. By abstract nonsense (\cite[Section2]{Bass}), since $F$ is (right) exact, the functor $F$ is isomorphic to
tensoring with the $A^\p$-bimodule $F(P^\p)$, however we claim the following
stronger result describing the natural transformations $\HOM(F,G)$
between projective functors $F$ and $G$:

\begin{theorem}
\label{Tcomb}
   Let $F$, $G$ be projective functors on $\MOF\text{-}A^\p$. Let
  $T\in\MOF\text{-}A^\p$ be a full projective tilting module. There is an isomorphism of
  vector spaces (even of rings if $F=G$)
  \begin{eqnarray*}
    \HOM(F,G)&\iso&\HOM_{D^\p\text{-}\MOF\text{-}A^\p}(F(T),G(T))\\
\phi&\mapsto&\phi_T.
  \end{eqnarray*}
\end{theorem}
\begin{proof}
  By naturality of the transformation $\phi$, the map $\phi_T$ is in fact a
  $D^\p\text{-}A^\p$-bimodule morphism and therefore, our map is
  well-defined. We first show that it is injective. Assume $\phi_T=0$. Let $P\in\MOF\text{-}A^\p$ be
  projective. There is an injective-tilting copresentation as in Lemma~\ref{copres}. Since $F$ and $G$ are
  exact and commute with finite direct sums we get a commutative diagram with
  exact rows of
  the form
  \begin{eqnarray*}
\xymatrix{
    0\ar@{->}[r]&F(P)\ar@{->}[r]\ar@{->}[d]^{\phi_P}&F(\bigoplus T)\ar@{->}[d]^0\\
    0\ar@{->}[r]&G(P)\ar@{->}[r]&G(\bigoplus T).
}
  \end{eqnarray*}
Since the vertical map on the right hand side is zero, we have $\phi_P=0$ as
well. Using a projective resolution we get $\phi_M=0$ for any
$M\in\pO$. Therefore, the map $\phi\mapsto\phi_T$ is injective.\\
Let $\phi_T\in\HOM_{D^\p\text{-}\MOF\text{-}A^\p}(F(T),G(T))$ and let $P\in\MOF\text{-}A^\p$ be projective with an
injective-tilting copresentation
\begin{eqnarray*}
  0\rightarrow P\stackrel{f}\lra\bigoplus_I T\stackrel{g}\lra\bigoplus_J T,
\end{eqnarray*}
Since $F$ and $G$ commute naturally with $\bigoplus$, we get a diagram of the
form
\begin{eqnarray*}
\xymatrix{
0\ar[r]&F(P)\ar@{->}[r]^{F(f)}&F(\bigoplus T)\ar@{->}[r]^{F(g)}\ar[d]^\wr&F(\bigoplus T)\ar@{->}[d]^\wr\\
&&\bigoplus F(T)\ar[d]^{\bigoplus\phi_T}&\bigoplus F(T)\ar[d]^{\bigoplus\phi_T}\\
&&\bigoplus G(T)\ar[d]^\wr&\bigoplus G(T)\ar@{->}[d]^\wr\\
0\ar[r]&G(P)\ar@{->}[r]^{G(f)}&G(\bigoplus T)\ar@{->}[r]^{G(g)}&G(\bigoplus T),}
\end{eqnarray*}
where the rows are exact. The isomorphisms exist and are natural,
since projective functors commute with direct sums (via the
natural isomorphisms $(\oplus M)\otimes E\cong \oplus (M\otimes
E)$ for $M\in\pO$ and $E$ a finite dimensional module). Since
$\phi_T$ is a $(D^\p,A^\p)$-bimodule morphism, it follows
that the rectangle in the diagram above commutes. Since the rows
are exact, restriction of the first vertical composition to $F(P)$
induces a unique morphism
$\phi_P\in\HOM_{\MOF\text{-}A^\p}(F(P),G(P))$. Standard arguments
show that $\phi_P$ does not depend on the chosen representation
and defines in fact a natural transformation $\phi$ between $F$
and $G$ restricted to the category of projective right
$A^\p$-modules. For $N\in\MOF\text{-}A^\p$ arbitrary we choose a
projective resolution $P^\bullet$. The naturality of $\phi$
defines a morphism of complexes $\Phi_N:F(P^\bullet)\lra
G(P^\bullet)$ inducing a unique map
$\phi_N\in\HOM_{\MOF\text{-}A^\p}(F(N),G(N))$ by exactness of $F$
and $G$. Again, standard arguments show that $\phi_N$ is
independent of the chosen projective resolution. Moreover, this
construction is natural in the sense that it defines a natural
transformation between $F$ and $G$ as functors on
$\MOF\text{-}A^\p$.
\end{proof}

Note that $\mV^\p F(T)$ is a $D^\p$-bimodule, where $df=F(d)\circ f$ and
 $fd=f\circ d$ for $d\in D^\p$, $f\in\mV^\p F(T)$. From
 Corollary~\ref{faithful} we directly get the following
\begin{cor}
\label{new} With the assumptions of Theorem~\ref{Tcomb} there is
an isomorphism of vector spaces
  \begin{eqnarray*}
    \HOM(F,G)&\iso&\HOM_{D^\p\text{-}\MOF\text{-}D^\p}(\mV^\p F(T),\mV^\p G(T))\\
\phi&\mapsto&\mV^\p(\phi_T).
  \end{eqnarray*}
\end{cor}

\subsection{The graded version}
To make the result stronger we would like to work in a graded setup. From
\cite{BGS}, it is known that $A^\p$ can be equipped with a non-negative
grading turning it
into a Koszul algebra. (In the example~\ref{ex} the grading is given by
putting all the arrows in degree $1$).

We call a module $\tilde{P}\in\gMOF\text{-}A^\p$ {\it a graded lift} of
$P\in\MOF\text{-}A^\p$ if it is isomorphic to $P$ after forgetting the
grading. Projective modules, standard objects and simple modules, for each of them exists a graded
lift (\cite{StGrad} or for a more general setup \cite{BinZhu}). Moreover, these lifts are unique up to isomorphism and
grading shift (\cite[Lemma 1.5]{StGrad}), since all these modules are
indecomposable. We fix {\it standard lifts} with the property that
their heads are concentrated in degree zero. In \cite{StGrad}, graded lifts of
translation functors are defined. By a graded lift we mean the following:

Assume
$C$, $B$ are graded rings and let $F:C\text{-}\MOD\rightarrow B$-$\MOD$ be a functor. Then a
functor $\tilde{F}: C\text{-}\gMOD\rightarrow B$-$\gMOD$ is a graded lift of $F$ if it is a $\mZ$-functor (i.e. it commutes with the
grading shifts in the sense of \cite[E.3]{AJS}),
such that $\op{f}_B\tilde F\cong F\op{f}_C$, where
$\op{f}_C:\gMOF\text{-}C\rightarrow \MOF\text{-}C$, $\op{f}_B:\gMOF\text{-}B\rightarrow\MOF\text{-}B$
denote the functors which forget the grading.

In \cite{StGrad}, it is shown that the translation functors
$\theta_s:\MOF\text{-}A^\p\rightarrow \MOF\text{-}A^\p$, for any simple reflection s, have graded lifts
\begin{eqnarray*}
\tilde\theta_s:\gMOF\text{-}A^\mb\rightarrow\gMOF A^\mb.
\end{eqnarray*}
On the other hand the natural projection of Koszul algebras
$A\rightarrow A^\p$ is graded and therefore, the functor $\tilde\theta_s$
restricts to a graded lift of $\theta_s:\MOF\text{-}A^\p\rightarrow\MOF\text{-}A^\p$ for any
$\p$. The following statement can be obtained as a direct consequence of the classification theorem of projective
functors (\cite{BGS}). We give an easy direct proof. Recall that a functor $F$
between abelian categories
is {\it indecomposable} if $F\cong F_1\oplus F_2$ implies $F_i=0$ for at least
one $i\in\{1,2\}$.

\begin{lemma}
\label{indec}
  For any simple reflection $s$, the functor
  $\theta_s:\MOF\text{-}A^\mb\rightarrow\MOF A^\mb$ is indecomposable.
  A graded lift $\tilde\theta_s$ is unique up to isomorphism and
grading shift.
\end{lemma}
\begin{proof}
Let $F=F_1\oplus F_2$ such that $F_i\not=0$ for $i=1,2$. In
particular, $FM=F_1M\oplus F_2M$ for any standard module $M$. On
the other hand, $\theta_s\Delta(x\cdot0)\in\cO_0$ is
indecomposable. (From the properties of translation functors it
follows easily that the socle of $\theta_s(M)$ is simple, hence
$\theta_s M$ is indecomposable.) In particular, $F_{i(M)}(M)=0$
for some $i(M)\in\{1,2\}$. The socle of any standard module
$\Delta(x\cdot0)\in\cO^\mb$ is of the form $\Delta(w_0\cdot0)$,
hence not annihilated by any $\theta_s$. From the exactness of $F$
we get that $i:=i(\Delta(w_0\cdot0))=i(M)$ for any standard module
$M$. Hence $F_{i}$ is zero when restricted to the category of
modules with standard flags. Since any projective module has a
standard flag, the functor is zero on projectives, hence vanishes
completely, since it is exact. Therefore, $F_{i}=0$. This
contradicts our assumption and therefore $F$ is indecomposable.
Since $F$ is (right) exact, it is given by tensoring with some
$A^\mb$-bimodule $X$. Since $X$ is indecomposable, a graded lift
is unique up to isomorphism and shift in the grading (this is
\cite[Lemma 1.5]{StGrad} applied to $A^\mb\otimes (A^\mb)^{opp}$).
\end{proof}

\begin{remark}
{\rm
Lemma~\ref{indec} only holds in the case $\p=\mb$. In general, the restriction of
  an indecomposable projective functor on $\cO_0^\mb$ to a functor on $\cO_0^\p$ could be
  decomposable or even zero (see \cite[Examples 3.7, Theorem 5.1]{StDuke}). In general it is not known
  how the restrictions decompose, not to mention a possible
  classification.}
\end{remark}

For any simple reflection $s$ we fix a standard lift $\tilde\theta_s$ of $\theta_s:\gMOF\text{-}A^\p\rightarrow\gMOF\text{-}A^\p$ such that
the standard lift of
$\Delta(0)$ is mapped to the standard lift of $P(s\cdot0)$. An endofunctor of
$\gMOF\text{-}A^\p$ is called {\it graded projective},
if it is a direct sum of (graded) direct summands of compositions of graded
lifts of translation functors. If $M$, $N\in\gMOF\text{-}B$ for some graded ring $B$, then $\HOM_B(M,N)$ is graded
by putting $$\HOM_B(M,N)_n=\{f\in\HOM_{B}(M,N)\mid f(M_k)\subseteq M_{k+n},
\forall k\in\mZ\}.$$ If $C$ is also a graded ring and
$F:\gMOF\text{-}B\rightarrow\gMOF\text{-}C$ is an exact functor then $F$ is given by
tensoring with some graded $B-C$-bimodule $X_F$. If $G$ is another such functor we
set
\begin{eqnarray*}
  \HOM(F,G)&=&\bigoplus_{n\in\mZ}\HOM(F,Q)_n=\bigoplus_{n\in\mZ}\HOM_{B\otimes
  C^{\op{opp}}\text{-}\MOF}(X_F,X_G)_n.
\end{eqnarray*}

We get the following refinement of
Theorem~\ref{Tcomb}:
\begin{theorem}
\label{Tgradcomb}
Let $F$, $G:\gMOF\text{-}A^\p\lra\gMOF\text{-}A^\p$ be graded projective functors. There is
an isomorphism of graded vector spaces (even of rings if $F=G$)
\begin{eqnarray*}
  \HOM(F,G)&\cong&\HOM_{D^\p\text{-}\MOF\text{-}A^\p}(F(\tilde{T}), G(\tilde{T}))\\
\phi&\mapsto&\phi_{\tilde{T}}.
\end{eqnarray*}
\end{theorem}
\begin{proof}
  Let $\Delta^\p(0)\in\gMOF\text{-}A^\p$ be the standard lift of the projective standard object. There is an inclusion $\Delta^\p(0)$ into
  $\bigoplus_{i\in I}\tilde{T}\langle i\rangle$ for a finite multiset $I$
  with elements from $\mZ$, because the injective hull of $\Delta^\p(0)$ is a
  direct sum of indecomposable projective-injective modules (compare the proof
  of Lemma~\ref{copres}). Since $\tilde{T}$ has a graded Verma flag
  (\cite[Theorem 7.2]{StGrad}), the
  cokernel of this inclusion has again a graded Verma flag by weight
  considerations. Therefore $P=\Delta^\p(0)$ has an injective-projective
  resolution of graded right $A^\p$-modules, i.e.\ there is an exact sequence of
  graded modules of the form
  \begin{eqnarray}
\label{gradedcopres}
     0\rightarrow P\lra\bigoplus_{i\in I}\tilde{T}\langle
     i\rangle\lra\bigoplus_{j\in J}\tilde{T}\langle j\rangle
  \end{eqnarray}
for some finite multisets $I$ and $J$ with entries in $\mZ$. The
statement follows then analogously to Proposition~\ref{copres} and
Theorem~\ref{Tcomb}.
\end{proof}

We denote by $\cZ(R)$ the centre of any ring $R$. We want to give
at least some description of the centre of
$\pO\cong\gMOF\text{-}A^\p$ which is by definition the centre of
the ring $A^\p$. Note that it inherits a grading from $A^\p$.
\begin{cor}
\label{center}
  Let $\ID$ denote the identity functor on $\gMOF\text{-}A^\p$. There are
  isomorphisms of (graded) rings
  \begin{eqnarray*}
    \END(\ID)\stackrel{(1)}\cong\END_{D^\p\text{-}\MOF\text{-}A^\p}(\tilde{T})\stackrel{(2)}\cong\cZ(D)\stackrel{(3)}\cong\cZ(A^\p).
  \end{eqnarray*}
In particular, the homogeneous part of degree zero in $\cZ(A^\p)$ has
dimension one.
\end{cor}
\begin{proof}
  The existence of the first two isomorphism follows directly from the previous
  two theorems and Corollary~\ref{new}. The isomorphism
 (3) is just obtained from the natural isomorphism $\END(\ID)\cong\cZ(A^\p)$ given by
  $\phi\mapsto\phi_{A^\p}(1)$. The last statement follows directly from the
  definition of the grading on $A^\p$, since $A^\p$ is indecomposable and
  its homogeneous part of degree zero is semisimple.
\end{proof}

\label{x}
Let $x\in W$ and $[x]=s_1s_2\cdot\ldots\cdot s_r$ be a fixed composition of
simple reflections. We denote by
$\theta_{[x]}=\theta_{s_r}\cdot\ldots\cdot\theta_{s_2}\theta_{s_1}$
the corresponding composition of translation functors and its
graded version
$\tilde\theta_{[x]}=\tilde\theta_{s_r}\cdots\tilde\theta_{s_2}\tilde\theta_{s_1}$,
where $\tilde\theta_s$ denotes the standard graded lift of
$\theta_s$. Let $R$ be a $\mC$-algebra with a $W$-action. For
$s\in W$ a simple reflection let $R^s$ be the invariants under
$s$. For $x\in W$ as above we denote
$R_{[x]}=\bullet\otimes_{R^{s_1}}R\otimes_{R^{s_{2}}}\otimes\cdots\otimes_{R^{s_r}}R$
considered as a functor on $\MOD\text{-}R$. If additionally $R$ is
graded and the action of $W$ is homogeneous, then $R_{[x]}$ is
even an endofunctor of $\gMOD\text{-}R$. In case $\cC$ is any
category having direct sums and $\cA$ is a list of objects of
$\cC$ then we write $\op{Ind}_\cC(\cA)$ for the set of iso-classes
of direct
summands of elements in $\cA$.

The next theorem is the classification
theorem from~\cite{BG}. We indicate a proof using Theorem~\ref{Tcomb},
Corollary~\ref{new} and the deformation theory from \cite{SHC}:

\begin{theorem}
\label{class}
  Let $\p=\mb$. There are natural bijections of isomorphism classes
\begin{eqnarray*}
  \begin{array}[t]{cc}
    \big\{\text{indecomposable projective functors on} \MOF\text{-}A^\mb\}&F\\
\updownarrow \text{\small 1:1}&\downarrow\\
\op{Ind}_{D^\mb\text{-}\MOF\text{-}A^\mb}(\theta_{[x]}T, x\in W) &F(T)\\
\updownarrow \text{\small 1:1}&\downarrow\\
\big\{\text{indecomposable projective objects of }\MOF\text{-}A^\mb\big\}
&F(\Delta^\mb(0))
  \end{array}
\end{eqnarray*}
\end{theorem}
\begin{proof}
 Theorem~\ref{Tcomb} implies that there is a
bijection between indecomposable projective functors on $\MOF\text{-}A^\mb$ and the set $\op{Ind}_{D^\mb\text{-}\MOF\text{-}A^\p}(\theta_{[x]}T, x\in W)$, where $T$ is the
indecomposable projective-injective module in $\MOF\text{-}A^\mb$. Consider
the algebra $S=\cU(\mh)$ as a graded algebra with $S_2=\mh$. Set
$S_+=(\mh)$. By Soergel's Endomorphismensatz (\cite{Sperv}) we have a
canonical isomorphism $D^\mb\cong S/(S_+^W)$, the coinvariant algebra which we
denote by $C$. By Corollary~\ref{new}, the functor $\mV^\mb$ defines a bijection between
the isomorphism classes of indecomposable functors on $\MOD{\text-}A^\p$ and
$$\op{Ind}_{C\text{-}\MOF\text{-}D^\mb}(\mV^\mb\theta_{[x]}T, x\in
W)=\op{Ind}_{C\text{-}\MOF\text{-}C}(C_{[x]}(C),x\in W),$$
since it is known from \cite{Sperv} that for any simple reflection $s$, there is an isomorphism of functors
\begin{eqnarray}
  \label{TV}
\mV^\mb\theta_s\cong \mV^\mb(\bullet\otimes_{C^s}C).
\end{eqnarray}
There are isomorphisms of $C$-bimodules, or $S$-bimodules as
follows:
$$C_{[x]}(C)=\big(S/(S_+^W)\big)_{[x]}(S/(S_+^W) =S_{[x]}(S/(S_+^W),$$ since
$S^W\subseteq S^s$  for any simple reflection $s$, and
$S_{[x]}(S)(S/(S_+^W)=S_{[x]}(S)\otimes_S (S/(S_+^W)$. Now, the
$S$-bimodules $S_{[x]}(S)$ are exactly the {\it Soergel special
  bimodules} as defined in \cite{SHC} and \cite[Bemerkung 5.12]{Sbimods}. By \cite[Proposition 11]{SHC}
specialisation gives a bijection between
$\op{Ind}_{S\text{-}\MOF\text{-}S}(S_{[x]}(S),x\in W)$ and
$\op{Ind}_{C\text{-}\MOF\text{-}C}(C_{[x]}(C),x\in W)$. The
theorem follows therefore directly from \cite[Section 6]{Sbimods}.
\end{proof}

\begin{rk}{\rm
  \begin{enumerate}[a.)]
  \item Let $\p=\mb$. Since all the constructions are compatible with gradings
  we get also a natural bijection of iso-classes of indecomposable graded projective endofunctors
  and indecomposable projective objects of $\gMOF\text{-}A$.
\item Theorem~\ref{class} is not true for arbitrary $\p$ (see \cite[Examples 3.7]{StDuke}).
\item Using the results from Kazhdan-Lusztig theory one can view
  Theorem~\ref{class} as a categorification of the integral group algebra
  $\mZ[W]$ as follows: There is an isomorphism of
  $\mZ$-algebras between $\mZ[W]$ and the Grothendieck ring of the
  indecomposable projective functors on $\cO_0^\mb$ mapping the Kazhdan-Lusztig
  basis element $C_w$ in the notation of \cite{SoKipp} to the indecomposable
  projective functor which maps the projective Verma module in $\cO_0^\mb$ to
  the indecomposable projective module with highest weight $w\cdot0$. More
  generally, the graded versions of projective functors categorify the
  (generic) Hecke algebra corresponding to the Weyl group of $\mg$ (a precise
  formulation can be found in \cite[Corollary 2.5]{StDuke}).
\end{enumerate}
}
\end{rk}
\section{An invariant of tangle cobordisms}
\label{inv} Now we would like to go one step further in the
categorification program proposed in \cite{BFK}. The main result
of \cite{StDuke} describes a functorial tangle invariant in terms
of bounded derived categories of certain $\cO_0^\p$ for $\mg$ of
type $A$. More precisely we associated to each $(m,n)$-tangle
diagram a functor between a bounded derived category defined by
certain $\cO_0^\p$ for $\mathfrak{sl}_m$ and a bounded derived
category defined by certain $\cO_0^\p$ for $\mathfrak{sl}_n$. It
was shown that, up to shifts, this assignment defines a functorial
invariant. In this section we will prove what was already
announced in \cite{StDuke}, namely that this assignment can be
extended firstly to a functorial invariant of oriented tangles
(such that the discrepancy with respect to shifts disappears) and
secondly to a $2$-functor, which means one can associate to each
oriented cobordism between two oriented tangles a (up to scalars)
well-defined homomorphism between the corresponding functors. In
this way the functorial tangle invariant is extended to an
invariant of tangles with cobordisms.

\subsection{The relevant functors and their morphisms}
Let now be $\mg=\SL_{n}$. We have $W=S_{n}$ generated by the simple
reflections $s_i$, $1\leq i<
n$ with the relations $s_is_j=s_js_i$ if $|i-j|\geq 2$ and
$s_is_js_i=s_js_is_j$ if $|i-j|=1$. For $1\leq k <n$ let $\p_k$ denote
the maximal parabolic subalgebra corresponding to the simple root
$\alpha_k$ (i.e. the corresponding parabolic subgroup is generated by all
$s_j$, where $j\not=k$). Set $\p_0=\p_n=\mg$. To simplify notation let
$\cO^{\p_i}_\mu$ denote the zero category if $n<i$ or $i<0$ and $\mu\in\mh^\ast$.\\

Fix $1\leq i<n$, $k\in\{0,\ldots, n\}$. Let $\theta_i=\theta_{s_i}$
denote the translation functor through the $i$-th wall with its
standard lift $\tilde\theta_i$. Let $\la_i\in\mh^\ast$ be integral
with stabiliser $W_{\la_i}=\{e,s_i\}$. In the following we need
the translation functors
\begin{eqnarray}
  \label{eq:outon}
  \theta_0^i:&&\cO_0^{\p_k}(\SL_n)\lra\cO_{\lambda_i}^{\p_k}(\SL_n)\\
  \theta^0_i:&&\cO_{\la_i}^{\p_k}(\SL_n)\lra\cO_0^{\p_k}(\SL_n)
\end{eqnarray}
{\it on} and {\it out off} the $i$-th wall. These are adjoint functors such
 that $\theta_i^0\theta_0^i\cong\theta_{i}$ (see e.g. \cite{GJ}). Enright and Shelton (\cite[chapter 11]{ES}) defined an equivalence of categories
\begin{eqnarray}
\zeta_{n,k}:\cO_{\la_1}^{\p_k}(\SL_n)\iso\cO^{\p_{k-1}}_0(\SL_{n-2}).
  \label{eq:ES}
\end{eqnarray}
We consider the following functors
\begin{eqnarray*}
\cap_{i,n}^{k}:\cO^{\p_k}_0(\SL_n)&\longrightarrow&\cO^{\p_{k-1}}_0(\SL_{n-2})\\
\cup_{i,n}^{k}:\cO^{\p_k}_0(\SL_n)&\longrightarrow&\cO^{\p_{k+1}}_0(\SL_{n+2})
\end{eqnarray*}
defined as
\begin{eqnarray*}
  \cap_{i,n}^{k}&=&\zeta_{n,k}\theta_0^1\theta_2\theta_3\cdots\theta_i\\
  \cup_{i,n}^{k}&=&\theta_i\theta_{i-1}\cdots\theta_2\theta_1^0\zeta^{-1}_{n+2,k+1}
\end{eqnarray*}
For $P^{\p_k}_n\in\cO^{\p_k}(\SL_n)$ a minimal projective generator we denote its endomorphism
ring by $A_n^k$. For the following we fix equivalences~\eqref{equiv} for any
$A_n^k$ and, concerning the notation, we will not distinguish the objects and
functors on each side. We have graded versions $\tilde\theta_i$,
$\tilde\theta_0^{\la_i}$, $\tilde\theta_{\la_i}^0$ of translation
through, on and out of the wall as defined in \cite{StGrad}. They are
normalised such that $\tilde\theta_0^{\la_i}$ (resp. $\tilde\theta_{\la_i}^0$) maps a simple module concentrated in degree zero to a
simple module (resp. a module with head) concentrated in degree $-1$.
Since the equivalences
$\zeta_{k,n}$ are compatible with the grading (\cite{Steen}), we get also (standard)
graded lifts $\tilde\cap_{i,n}^k$ and $\tilde\cup_{i,n}^k$ of $\cap_{i,n}^k$
and $\cup_{i,n}^k$ respectively. They have the following properties:
\begin{lemma}
\label{capcup}
  \begin{enumerate}[a.)]
  \item The functors $\cap_{i,n}^k$
and $\cup_{i,n}^k$ are indecomposable. In particular, a graded lift is unique
up to isomorphism and grading shift.
\item The standard lifts provide adjoint pairs of functors
  \begin{eqnarray*}
    \big(\tilde\cap_{i,n}^k\langle -1\rangle,\tilde\cup_{i,n-2}^{k-1}\big)\quad\text{and}\quad
\big(\tilde\cup_{i,n-2}^k\langle 1\rangle,\tilde\cap_{i,n}^{k+1}\big).
  \end{eqnarray*}
  \end{enumerate}
\end{lemma}
\begin{proof}
  \begin{enumerate}[a.)]
  \item Assume that $\cap_{i,n}^k$ decomposes as $\cap_{i,n}^k\cong F_1\oplus F_2$. Then
  $\cup_{i,n-2}^{k-1}\cap_{i,n}^k\cong\cup_{i,n-2}^{k-1}F_1\oplus\cup_{i,n-2}^{k-1}F_2$.
  Since $\cup_{i,n-2}^{k-1}\cap_{i,n}^k\cong\theta_i$ is indecomposable
  (\cite[Lemma 6.3, Theorem 5.1]{StDuke}) it follows say
  $\cup_{i,n-2}^{k-1}F_1=0$. Assume $\cup_{i,n-2}^{k-1}(L)=0$ for some simple
  object $L$. Then
  $0=\theta_{i-1}\cup_{i,n-2}^{k-1}(L)\cong\theta_{i-1}\theta_i\cdots\theta_2\theta_1^0(L')$
  for some simple object $L'$. Since
  $\theta_j\theta_{j'}\theta_j\cong\theta_{j}$ for $|j-j'|=1$ (\cite[Theorem 4.1]{StDuke}) it follows
  inductively that $\theta_2\theta_1^0(L')=0$, hence
  $0=\theta_0^1\theta_2\theta_1^0(L')\cong L'$ because of
  \cite[Lemma 4]{BFK}. Therefore, $F_1=0$ contradicting our assumption. That
  means $\cap_{i,n}^{k}$ is indecomposable. The uniqueness of a graded lift follows form
  \cite[Lemma 1.5]{StGrad}
  (applied to the rings $A_n^k\otimes {(A_{n-2}^{k-1})}^{\op{opp}}$ and $A_n^k\otimes
  {(A_{n+2}^{k+1})}^{\op{opp}}$ respectively).\\
Assume $\cup_{i,n}^{k}$ decomposes as $\cup_{i,n}^k\cong F_1\oplus F_2$ is decomposable. As above we
  deduce that, without loss of generality, $F_1\cap_{i,n+2}^{k+1}=0$. Hence,
  $0=F_1\cap_{i,n+2}^{k+1}\cup_{i,n}^k\cong F_1\oplus F_1$ by \cite[Theorem
  6.2]{StDuke}. We get $F_1=0$ contradicting our assumption.
\item The adjointnesses follow directly from the definitions and the
  adjoint pairs $(\tilde\theta_1^0\langle 1\rangle,\tilde\theta_0^1)$,
  $(\tilde\theta_0^1\langle -1\rangle,\tilde\theta_1^0)$, and
  $(\tilde\theta_i,\tilde\theta_i)$ (\cite[Theorem 8.4, Corollary 8.3]{StGrad}).
\end{enumerate}
\end{proof}

We describe the endomorphism rings of these functors.

\begin{theorem}
\label{Endcapcup}
For any $1\leq i\leq n$ the following holds:
  \begin{enumerate}[a.)]
  \item \label{a}
There is an isomorphism of graded vector spaces
\begin{equation*}
\END(\tilde\cup_{i,n}^k)\cong\cZ(A_n^k)\oplus\cZ(A_n^k)\langle
  2\rangle\cong\END(\tilde\cap_{i,n}^k).
\end{equation*}
Hence they are non-negatively
  graded and one-dimensional in degree zero.
\item  \label{b} The vector spaces $\HOM(\ID, \tilde\theta_i)$ and $\HOM(\tilde\theta_i,\ID)$ are both strictly positively
  graded and one-dimensional in degree one (with basis the adjunction morphism).
  \end{enumerate}
\end{theorem}
\begin{proof}
There are isomorphisms of graded vector spaces
    \begin{eqnarray*}
      \begin{array}[t]{lclll}
     \HOM(\tilde\cup_{i,n}^k,\tilde\cup_{i,n}^k)&\cong&\HOM(\ID\langle-1\rangle,\tilde\cap_{i,n+2}^{k+1}\tilde\cup_{i,n}^k)&\text{(Lemma~\ref{capcup})}\\
&\cong&\HOM(\ID\langle -1\rangle,\ID\langle -1\rangle\oplus\ID\langle 1\rangle)
&\text{(\cite[Theorem 6.2]{StDuke})}\\
&\cong&\HOM(\ID,\ID\oplus\ID\langle 2\rangle).
   \end{array}
\end{eqnarray*}
The existence of the first isomorphism follows therefore from
Corollary~\ref{center}. Assuming the existence of the second
isomorphism from part~\eqref{a}, part~\eqref{b} follows by
adjunction (Lemma~\ref{capcup}) and the self-adjointness of
$\tilde\theta_i$ (\cite[Corollary 6.3]{StGrad}), since
$\tilde\theta_i\cong\cup_{i,n-2}^{k-1}\tilde\cap_{i,n}^k$ (see
\cite[Proposition 6.7]{StDuke}), and the adjunction morphisms are
both of degree one (\cite{StGrad}). To establish the second
isomorphism of part~\eqref{a} we have to work more. Let for the
moment $P$ be the chosen minimal projective generator of
$\cO_0^{\p_k}(\SL_n)$ with endomorphism ring $A=A_n^k$ and let
$P_\la$ be a minimal projective generator of $\cO_\la^{\p_k}$
where $\la$ is an integral weight with stabiliser
$W_\la=\{e,s_1\}$. Set $B=\END_\mg(P_\la)$. The composition
$F=\theta_0^1\theta_2\cdots\theta_i:\MOF-A\rightarrow \MOF-B$ is
therefore given as tensoring with the bimodule
\begin{eqnarray*}
  {\bf F}:=\HOM_\mg(P_\la,F(P))\in A\text{-}\MOF-B
\end{eqnarray*}
with the actions $af=F(a)\circ f$ and $fb=f\circ b$, where $f\in{\bf F}$, $a\in
A$, $b\in B$.
Let $\hat{F}=\theta_i\cdots\theta_2\theta_1^0$ be the adjoint functor with
describing bimodule ${\bf \hat{F}}=\HOM_\mg(P,\hat{F}(P_\la))\in B\text{-}\MOF-A$.
Since $F\in A\text{-}\MOF$ and $B\in B\text{-}\MOF$, the space $\HOM_{\MOF-B}({\bf F},B)$
becomes naturally a $(B,A)$-bimodule. We claim that there are isomorphisms of $(B,A)$-bimodules
\begin{eqnarray*}
  \HOM_{\MOF-B}\big({\bf F}, B\big)\stackrel{\Phi}{\longleftarrow}\HOM_\mg\big(F(P),P_\la\big)\stackrel{\Psi}{\longrightarrow}{\bf\hat{F}},
\end{eqnarray*}
Here $(\Phi(g))(b)=g\circ b$ and $\Psi$ is given by adjointness, i.e. $g\mapsto \hat{F}(g)\circ\eta$, where
$\eta:\ID\rightarrow \hat{F}F$ is the adjunction morphism. We calculate
with $a\in A$ and $b\in B$ explicitly $\Psi(b g a)=\hat{F}(b g a)\circ\eta=\hat{F}(b g\circ
F(a))\circ\eta=\hat{F}(bg)\circ(\hat{F}F(a))\circ\eta=\hat{F}(bg)\circ\eta\circ
a=\hat{F}(b)\circ\Psi(g)\circ a=b\Psi(g)a$. Hence, $\Psi$ is in fact an
isomorphism of bimodules. Since $\Phi(bga)(h)=\Phi(b\circ g\circ
F(a))(h)=\Phi(b\circ g)(F(a)\circ h)=(b\Phi(g))(F(a)\circ h)=(b\Phi(g) a)(h)$,
the map $\Phi$ is compatible with the bimodule structures. It is obviously an
isomorphism, since $P_\la$ is a projective generator.\\
Now, $\HOM_{\MOF-B}(\bullet,B)$ defines an equivalence $e:A\text{-}\MOF-B\iso B\text{-}\MOF-A$
with inverse functor $\HOM_{B\text{-}\MOF}(\bullet, B)$. Therefore,
\begin{eqnarray*}
\END(F)\cong\END_{A\text{-}\MOF-B}({\bf
  F})\cong\END_{B\text{-}\MOF-A}({\bf\hat{F}})\cong\END(\hat{F}).
\end{eqnarray*}
From the definitions of the graded translation functors it follows directly
that all the isomorphisms are grading preserving. The Theorem follows.
\end{proof}

To simplify the setup, instead of working with derived categories, we will
work with homotopy categories of complexes. To make it consistent with
\cite{StDuke}, one only has to replace all bounded derived categories
appearing there by the bounded homotopy category of perfect complexes (for the
general setup we refer for example to \cite{KZ}). Given a complex
$(X^\bullet,d)$ we have the differentials $d:X^i\rightarrow X^{i-1}$ and
$(X[k])^i=X^{i+k}$ for any $k\in\mZ$.\\

Let $\fD^b_{per}(\gMOF-A_n^k)$ denote the bounded homotopy
category of perfect complexes of graded $A_n^k$-modules. Let
$\mathtt{C}_i^k$ be the endofunctor of $\fD^b_{per}(\gMOF-A_n^k)$ given
by tensoring with the complex of graded $A_n^k$-bimodules
\begin{eqnarray}
  \label{eq:Ci}
  \cdots\rightarrow 0\rightarrow A_n^k\langle 1\rangle\rightarrow\tilde\theta_i
A_n^k\rightarrow 0\rightarrow\cdots
\end{eqnarray}
where the map is the adjunction morphism (see
Theorem~\ref{Endcapcup}) and $\tilde\theta_i A_n^k$ is
concentrated in position zero. The functor $\mathtt{C}_i^k$
defines an auto-equivalence of $\fD^b_{per}(\gMOF-A_n^k)$
(\cite[Section 7]{StDuke}). Let $\mathtt{K}_i^k$ be its inverse,
that is the functor given by tensoring with the complex
\begin{eqnarray}
  \label{eq:Ki}
\cdots\rightarrow
0\rightarrow\tilde\theta_i A_n^k\rightarrow A_n^k\langle -1\rangle\rightarrow
0\rightarrow\cdots
\end{eqnarray}
If $F$ is a finite composition of functors of the form
$\mathtt{C}_i^k$, $\mathtt{K}_i^k$ for some fixed $k$, then we
have (via \eqref{eq:Ci} and \eqref{eq:Ki}) the corresponding complex,
say $X(F)$ of graded $A_n^k$-bimodules. We consider $X(F)$ as an
object in $\fD^b(A_n^k\text{-}\gMOF\text{-}A_n^k)$, the
homotopy category of complexes of graded $A_n^k$-bimodules and
denote by $\END(F)$ its endomorphism ring.

\begin{lemma}
\label{CK} Let $F$ be a finite composition of functors of the form
$\mathtt{C}_i^k$, $\mathtt{K}_i^k$ for some fixed $k$.
\begin{enumerate}[a.)]
\item  There is an isomorphism of graded rings $\END(F)\cong \cZ({A_n^k})$.
\item If $G=\tilde\cap_{k,n}F$ and there exists a grading preserving
  isomorphism $f:G\cong H$ for some functor
$H:\fD^b_{per}(\gMOF-A_n^k)\rightarrow\fD^b_{per}(\gMOF-A_{n-2}^{k-1})$. Then $f$ is unique
up to a scalar.
\end{enumerate}
\end{lemma}

\begin{proof}
 The first statement follows directly from \cite{Rickardcentre}. For a nice
 direct argument we refer to \cite[Proposition 1]{Khovanovcob}.
 For the second statement we assume there is another isomorphism $f'$, then $f'^{-1}\circ f\in\END(G)$ is
  of degree zero and by Theorem~\ref{Endcapcup} a scalar multiple of the
  identity.
\end{proof}

We set $\mathtt{C}_i:=\oplus_{k=0}^n\mathtt{C}_i^k$, considered as
an auto-functor of $\bigoplus_{k=0}^n\fD^b_{per}(\gMOF-A_n^k)$ and let
$\mathtt{K}_i=\oplus_{k=0}^n\mathtt{K}_i^k$ be its inverse.

\subsection{The tangle $2$-category and its generators}
Let $A_n=\oplus_{k=1}^n A^k_n$. In \cite{StDuke}, we assigned to a plane
  diagram of a tangle with $m$
  bottom and $n$ top points a functor from $\fD^b_{per}(\gMOF-A_m)$ to
  $\fD^b_{per}(\gMOF-A_n)$. Up to shifts, this assignment provided an invariant of
  isomorphism classes of tangles. We also looked at the category $\mathcal{COB}$ of
  $2$-cobordisms, where the objects are a finite number of labelled oriented
  one-manifolds and the morphisms are cobordisms. By considering the objects as
  special oriented $(0,0)$-tangles we assigned (in a $2$-functorial way) to each
  object some functor and to each morphisms a natural transformation
  between the corresponding functors (\cite[Theorem 8.1]{StDuke}) satisfying
  the defining relations for the isotopy classes of cobordisms. In the
  following we will show that this can be extended to arbitrary oriented
  tangles giving rise to a functorial invariant of tangles and cobordism.

More precisely, let ${\mathcal Tan}$ denote the category of tangles,
  i.e. objects are the positive integers and morphisms are unframed
  tangle diagrams. Let ${\mathcal Tan}^{or}$ be the $2$-category
  of oriented tangles and cobordisms, i.e. objects are the positive integers,
  morphisms are unframed oriented tangles and $2$-morphisms are
  diagrams of tangle cobordisms. For details see for example \cite{BL},
  \cite{CRS}, \cite{CS}, \cite{Fischer}. This $2$-category is of interest,
  since the $1$-morphisms give rise to an algebraic description of tangles
  (and hence of knots and
  links), whereas the $2$-morphisms describe compact surfaces smoothly embedded
  in $\mR^4$.

\begin{figure}[htbp]
    \centering
    \includegraphics{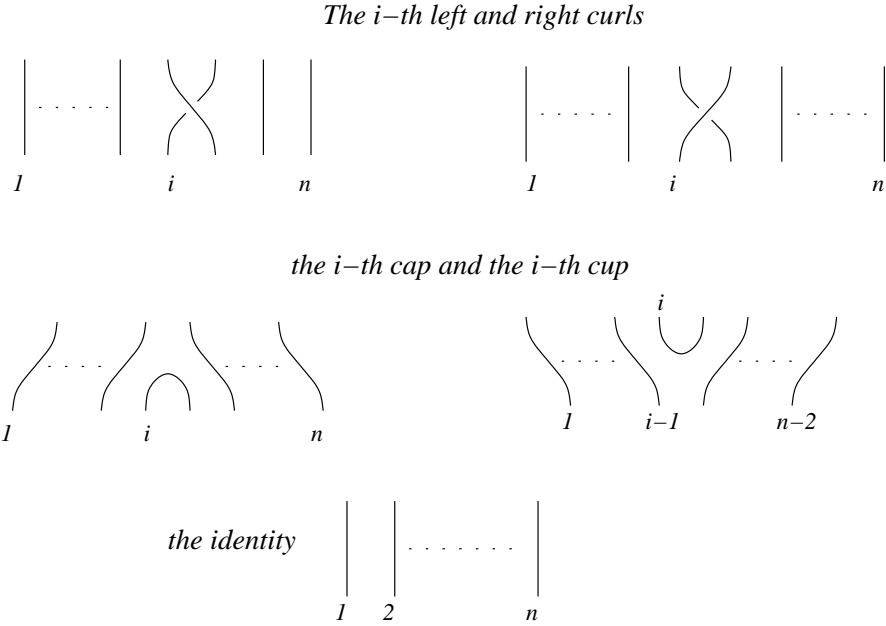}
    \caption{Elementary tangle diagrams generating the $1$-morphisms }
  \label{fig:elementary} 
 \end{figure}

\begin{figure}[htbp]
\begin{center}
    \includegraphics{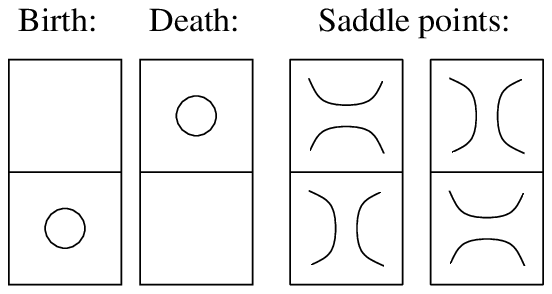}
    \includegraphics{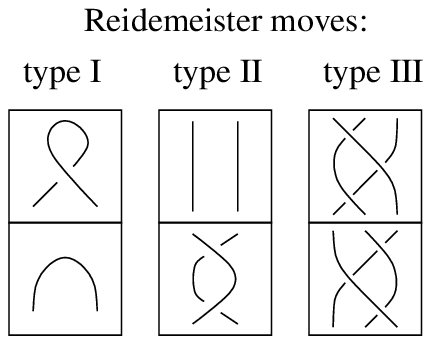}{\vspace {0.5cm}}
    \includegraphics{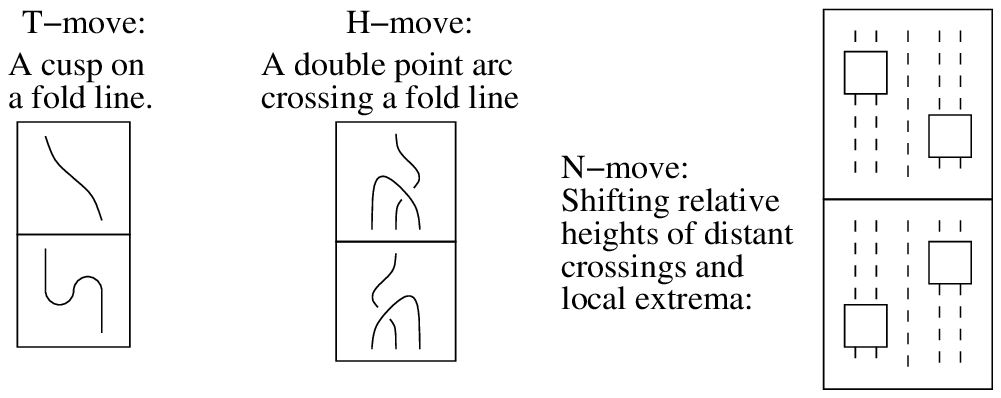}
\end{center}
\caption{Typical generating $2$-morphisms}
  \label{fig:2morphisms}
\end{figure}

The $1$-morphisms are generated by the {\it elementary tangles} as depicted in
  Figure~\ref{fig:elementary}.
i.e. the $1$-morphisms are just the products of elementary
  tangle diagrams. For an $(m,n)$-tangle diagram $T_1$ and an $(m',n')$-tangle diagram $T_2$, the composition $T_2T_1$ is defined if
  and only if $n=m'$. In this case $T_2T_1$ is an $(m,n')$-tangle diagram and obtained
  by putting $T_2$ on top of $T_1$ and identifying the bottom points of $T_2$
  with the top points of $T_1$. The $2$-morphisms in the category
  ${\mathcal Tan}$ are diagrams of cobordisms generated by birth, death, saddle points, Reidemeister
  moves, shifting relative heights of distant crossings, local extrema, the
  identity morphisms, cusps on fold lines, and double point arcs crossing a
  fold line. The {\it typical generators} (apart from the identities)
are depicted in Figure~\ref{fig:2morphisms}. These are the {\it
  elementary string interactions} from \cite{CRS}, where one can also find the
  corresponding surfaces displayed.
Any $2$-morphism is a composition of generating $2$-morphisms, the generating
$2$-morphisms are obtained by reading the typical generators either upwards or downwards, taking
their vertical and horizontal mirror images and changing between negative and
positive crossings. For details we again refer to \cite{CRS}.

\begin{figure}[htbp]
  \begin{center}
\includegraphics{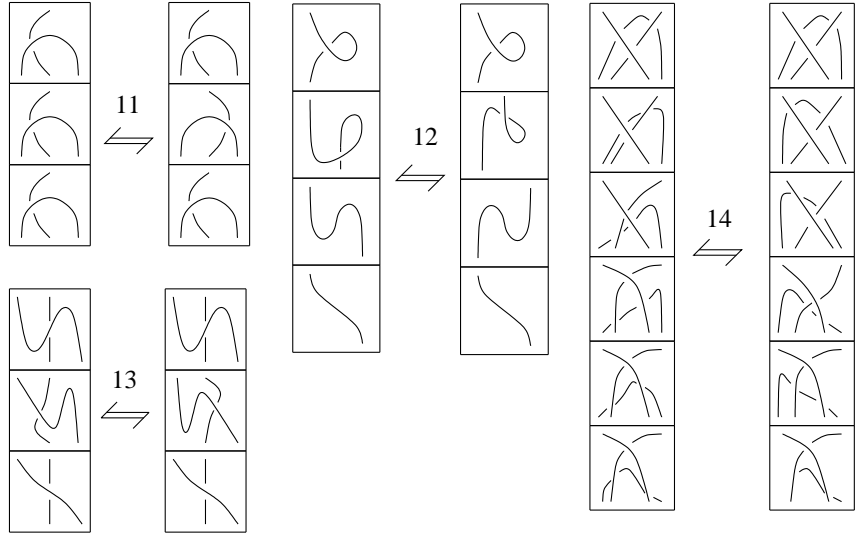}
\caption{The movie moves $11$ to $14$}
\label{fig:moves1114}
\end{center}
\end{figure}

Recall that two tangle diagrams represent (up to ambient isotopy) the same
tangle if they differ by a sequence of Reidemeister moves. As briefly mentioned in the introduction, to any cobordism (more precisely to any knotted surface
in the sense of \cite{CRS}) there is an associated sequence of generating $2$-morphisms. Moreover to any sequence of generating $2$-morphisms there is a
cobordism (knotted surface) whose diagram sequence is the given one
(\cite[Theorem 3.5.4]{CRS}).
 By \cite[Theorem 3.5.5]{CRS}, two sequences of generating $2$-morphisms represent (up to ambient
isotopy) the same cobordism if they differ by a sequence of so-called {\it
  movie moves} (see also \cite{CRS}, \cite{CS}, \cite{BL} or \cite{Fischer}). Examples of
movie moves are depicted in Figure~\ref{fig:moves1114}.
For a complete list of movie moves we refer to \cite[Theorem 3.5.5]{CRS}.

\subsection{The functorial invariant of tangles}
 As suggested in \cite{BFK}, we associated in \cite{StDuke} functors to
 elementary tangles as follows: To the identity tangle with $n$ strands we associate the identity functor on
$\fD^b_{per}(\gMOF\text{-}A_n)=\bigoplus_{k=0}^n\fD^b_{per}(\gMOF-A_n^k)$. For the U-turns we assign
\small
\begin{eqnarray*}
{\parbox{4.5cm}{\includegraphics{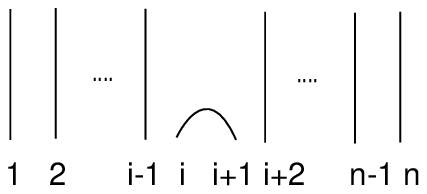}}}
\rightsquigarrow&&\tilde\cap_{i,n}\langle 1\rangle:\fD^b_{per}\big(\gMOF\text{-}A_n\big)\longrightarrow\fD^b_{per}\big(\gMOF\text{-}A_{n-2}\big).
\end{eqnarray*}
\begin{eqnarray*}
{\parbox{4.5cm}{\includegraphics{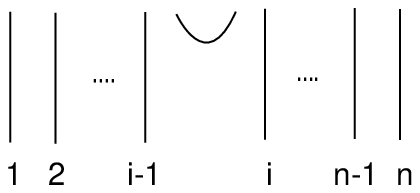}}}
\rightsquigarrow&&\tilde\cup_{i,n}\langle -1\rangle:\fD^b_{per}\big(\gMOF\text{-}A_n\big)\longrightarrow\fD^b_{per}\big(\gMOF\text{-}A_{n+2}\big).\\
\end{eqnarray*}
\normalsize
To the $i$-th right twisted curl
we associate the functor $\cC_i:=\mathtt{C}_i\langle 1\rangle$. To the $i$-th left twisted curl we
associate the inverse functor $\cK_i:=\mathtt{K}_i\langle -1\rangle$.

\begin{remark}
\label{renormalisation}
{\rm
  The assignments differ slightly from the ones in \cite{StDuke} as follows:
  In \cite{StDuke} we took $\tilde\cap_{i,n}$ instead of
  $\tilde\cap_{i,n}\langle 1\rangle$ and also $\tilde\cup_{i,n}$ instead of
  $\tilde\cap_{i,n}\langle -1\rangle$. Moreover, we swapped the assignments for
  left and right twisted curls. We introduced this renormalisation to make it
  compatible with Khovanov
 homology (see Conjecture~\ref{Khovlink}).
}
\end{remark}

Let ${\mathcal Fun}$ denote the category which we define as follows: The
  objects are the bounded homotopy categories $\fD^b_{per}(\gMOF\text{-}A_n)$
  (i.e. the objects are indexed by the natural numbers). The $1$-morphisms are
  functors between the corresponding categories.

\begin{theorem}(see \cite[Theorem 7.1]{StDuke})
\label{2cat}
  There is a functor
  \begin{eqnarray*}
    \Phi:{\mathcal Tan}\rightarrow{\mathcal Fun}
  \end{eqnarray*}
which is given on objects by
\begin{eqnarray*}
n\in\mN\mapsto{\fD}^b_{per}(\gMOF\text{-}A_n),
\end{eqnarray*}
and on elementary $1$-morphisms by the assignments above, such that
if $t_1$ and $t_2$ are $1$-morphisms which differ by a sequence of
  Reidemeister moves then there is an isomorphism of functors
  $\Phi(t_1)\cong\Phi(t_2)\langle 3r\rangle[r]$ for some $r\in\mZ$.
\end{theorem}

\begin{proof}
  This is \cite[Theorem 7.1]{StDuke} and
  Remark~\ref{renormalisation}, since the renormalisation is compatible with
  the Reidemeister, $N$, $H$, and $T$-moves.
\end{proof}

Hence, the functor $\Phi$ defines, up to shifts, a functorial
invariant of tangles. Moreover (see \cite[Proof of Theorem
7.1]{StDuke}), it turns out that the shift problems only occur in
$H$-moves.

\subsection{Oriented tangles and cobordisms}
Theorem~\ref{2cat} can be improved and made more natural by working with
{\it oriented} tangles and cobordisms instead. Let ${\mathcal Func}$ denote
the $2$-category with underlying category ${\mathcal Fun}$: The
  objects are the bounded homotopy categories $\fD^b_{per}(\gMOF\text{-}A_n)$, the $1$-morphisms are
  functors between the corresponding categories, the $2$-morphisms are the
  natural transformations between the functors, but after forgetting the
  grading and only up to a
  multiplication with a homogeneous element of degree $0$ of the centre of the
  source or image category. We would like to construct a functorial invariant
  of oriented tangles and cobordisms.

\begin{figure}[hbtp]
  \centering
  \includegraphics{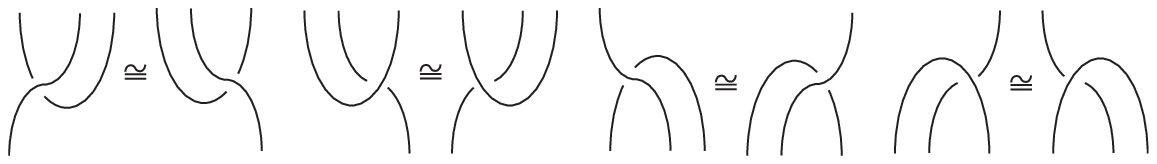}
  \caption{The $H$-moves}
  \label{fig:H-moves}
\end{figure}

To any elementary oriented tangle
without crossing we associate the same functors as before. Let us consider
the four $H$-moves of non-oriented tangles depicted in Figure~\ref{fig:H-moves},
and let $F_i$ $1\leq i\leq 8$ be the corresponding functors. Then it is known
 that we have isomorphisms of $\mZ$-functors $F_1\langle -3\rangle[-1]\cong F_2$,
$F_3\langle -3\rangle[-1]\cong F_4$, $F_5\langle 3\rangle[1]\cong F_6$,
$F_7\langle 3\rangle[1]\cong F_8$ (see \cite[Proof of Theorem 7.1]{StDuke},
 note the signs appear because of Remark~\ref{renormalisation}).

\begin{figure}[htbp]
\begin{tabular}{cc}
 \begin{minipage}{5cm}
   \centering
   \includegraphics{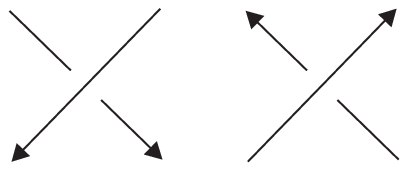}
   \caption{Oriented crossing inducing a shift $\langle -3\rangle[-1]$.}
   \label{fig:oriented2}
 \end{minipage} \hspace{1cm} &
 \begin{minipage}{5cm}
  \centering
  \includegraphics{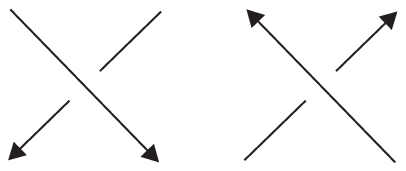}
  \caption{Oriented crossings inducing a shift $\langle3\rangle[1]$.}
    \label{fig:oriented1}
\end{minipage}
\end{tabular}
\end{figure}

For oriented crossings we modify the assignment as follows: To an oriented
crossing as depicted in Figure~\ref{fig:oriented2}

\begin{figure}[hbtp]
\centering
\includegraphics{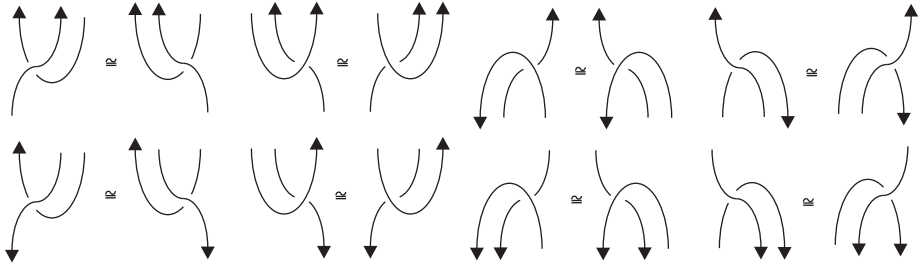}
 \caption{oriented $H$-moves}
\label{Hmovesor}
\end{figure}

we associate the functor $\cC_i\langle -3\rangle[-1]$ and to an oriented
crossing as depicted in Figure~\ref{fig:oriented1} we associate
the functor $\cK_i\langle 3\rangle[1]$. We get the oriented versions
of the $H$-moves as depicted in Figure~\ref{Hmovesor}. On can see immediately that our renormalisation precisely implies that the two
functors corresponding to a move are now isomorphic (without any shifts!).

\begin{theorem}
\label{2catorient}
  There is a functor of $2$-categories
  \begin{eqnarray*}
    \Phi^{or}:{\mathcal Tan}^{or}\rightarrow{\mathcal Func}
  \end{eqnarray*}
which is given on objects by
\begin{eqnarray*}
n\in\mN\mapsto{\fD}^b_{per}(\gMOF\text{-}A_n),
\end{eqnarray*}
and on elementary $1$-morphisms by the assignments above, such that
\begin{enumerate}
\item if $t_1$ and $t_2$ are $1$-morphisms which differ by a sequence of
  Reidemeister moves then there is an isomorphism of functors
  $\Phi^{or}(t_1)\cong\Phi^{or}(t_2)$.
\item if $c_1$ and $c_2$ are sequences of generating $2$-morphisms which
  differ by a sequence of movie moves then $\Phi^{or}(c_1)=\Phi^{or}(c_2)$.
\end{enumerate}
\end{theorem}

The proof of this result requires the following auxiliary lemma:

\begin{lemma}
\label{phi}
Let $n$ be any positive integer and let $1\leq i,k\leq n$.
There is an isomorphism of functors
$$\phi_{i,n}:\;\tilde\cap_{i,n+2}^{\p_k}\tilde\cup_{i,n}^{\p_k}\cong\ID\langle 1\rangle\oplus\ID\langle1\rangle:\gMOF-A_n^k\rightarrow\gMOF-A_n^k.$$
\end{lemma}

\begin{proof}
  Since $\p_k$ is a maximal parabolic, we have by \cite[Theorem 4.1]{StDuke}
  isomorphisms of endofunctors of $\gMOF-A_n^k$ as follows
  \begin{eqnarray}
\label{TL}
  \tilde\theta_i\tilde\theta_j\tilde\theta_i&\cong&\tilde\theta_j\tilde\theta_i\tilde\theta_j\quad\quad\quad\quad\text{
  if $|i-j|=1$},
\end{eqnarray}
It is well-known (see e.g. \cite{BG}) that
$\theta_0^i\theta_i^0\cong\ID\oplus\ID$ and $\tilde\theta_0^i\tilde\theta_i^0\cong\ID\langle
1\rangle\oplus\ID\langle-1\rangle$ (the latter by \cite[Theorem 8.2
(4)]{StGrad}). For $i=1$, the statement of the lemma follows then directly
from the definitions, since the Enright-Shelton equivalences are compatible
with the grading \cite{Steen}. For $i>1$ we can deduce, by applying the
relations \eqref{TL}, that it is enough to consider the case $i=2$. Then, by
\cite[Lemma 4]{BFK}, the lemma is true if we forget the grading. To determine
the graded shifts we just apply \cite[Theorem 8.2 (2), Theorem 5.3]{StGrad}
and the lemma follows.
\end{proof}

From now on we {\it fix} an isomorphism $\phi_{i,n}$ as in Lemma~\ref{phi} for any $i$, $n$.

\begin{proof}[Proof of Theorem~\ref{2catorient}]
From Theorem~\ref{2cat} we have the functor $\Phi$. Since
Reidemeister moves of type I do not involve any crossings of the
form displayed in Figure~\ref{fig:oriented1} or
Figure~\ref{fig:oriented2} we do not have to check anything there.
For type II moves we are either in the situation as for tangles
without orientation or we have the assignments from
Figure~\ref{fig:oriented1} and Figure~\ref{fig:oriented2} and
hence the statement follows. The moves of type III can be easily
checked, since we know the result for the non-oriented case. For
$T$-moves and $N$-moves nothing is to check. Finally, the functors
are defined such that the $H$-moves work. Therefore $\Phi^{or}$
defines a functorial invariant of oriented tangles.

We have now to assign to
each generating $2$-morphism a functor between the corresponding functors.

{\it Reidemeister moves:} All different types of Reidemeister moves correspond
to two isomorphic functors (\cite[Proof of Theorem
7.1]{StDuke}). For any such move $m$ we fix an isomorphism $i(m)$ which is
grading preserving such that if we read the move upside down we get the
inverse isomorphism. Then we define
$\Phi(m)=i(m)$.

{\it T-moves:} By \cite[Theorem 6.2, (6.1), (6.2)]{StDuke} we know that all
the T-moves correspond to pairs of isomorphic functors if we forget the
grading. Assume
$\tilde\cap_{i+1,n+2}^{\p_k}\tilde\cup_{i,n}^{\p_k}\cong\ID\langle k\rangle$
for some $k\in\mZ$. Then we have
$\tilde\cup_{i+1,n}^{\p_k}\tilde\cap_{i+1,n+2}^{\p_k}\tilde\cup_{i,n}^{\p_k}\tilde\cap_{i,n+2}^{\p_k}\cong\tilde\cup_{i+1,n}^{\p_k}\tilde\cap_{i,n+2}^{\p_k}\langle
k\rangle$. However \cite[Proposition 6.7]{StDuke} tells us that the left hand
side of this isomorphism is isomorphic to
$\tilde\theta_{i+1}\tilde\theta_{i}$, whereas the right hand side is
isomorphic to $\tilde\theta_{i+1}\tilde\theta_{i}\langle k\rangle$. Hence
$k=0$. For any typical T-move $m$ we fix a grading preserving isomorphism
$t(m)$ such that if we read the move upside down we have the inverse
isomorphism. Then we define $\Phi(m)=t(m)$.

{\it H-move:} All different types of $H$-moves correspond
to two isomorphic functors. For any such move $m$ we fix an isomorphism $h(m)$ which is
grading preserving such that if we read the move upside down we get the
inverse isomorphism. Then we define
$\Phi(m)=h(m)$.

{\it N-moves:} All different types of $N$-moves correspond
to two isomorphic functors. For any such move $m$ we fix an isomorphism $n(m)$ which is
grading preserving such that if we read the move upside down we get the
inverse isomorphism. Then we define
$\Phi(m)=n(m)$.

{\it Saddle point:} There are the degree preserving adjunction morphisms
$\ID\langle 1\rangle\rightarrow\tilde\theta_i$ and
$\tilde\theta_i\rightarrow\ID\langle -1\rangle$
(Theorem~\ref{Endcapcup}). We fix isomorphisms
$\psi_{i,n}:\tilde\cup_{i,n}\tilde\cap_{i,n}\cong\tilde\theta_i$ which exist by \cite[Proposition
6.7]{StDuke}. To each saddle point move we associate the natural
transformation which is induced from the corresponding adjunction morphism.

{\it Births/Deaths:} Via the isomorphisms $\phi_{i,n}$ we get a surjection
$\tilde\cap_{i,n+2}^{\p_k}\tilde\cup_{i,n}^{\p_k}\rightarrow \ID\langle
1\rangle$ and an inclusion $\ID\langle -1\rangle\rightarrow
\tilde\cap_{i,n+2}^{\p_k}\tilde\cup_{i,n}^{\p_k}$ (both maps homogeneous of degree zero). To each birth move we
associate the corresponding natural transformation (which is then homogeneous
of degree $-1$).

It is left to show that the natural transformations satisfy the relations
given by the movie moves. Because of our results in Section~\ref{inv}, the arguments are quite routine and mimic the
arguments in \cite{Khovanovcob}. Now the statement follows directly by copying the arguments from \cite{Khovanovcob} if we make the following correspondences:
Proposition 2 there corresponds to our Corollary~\ref{center}; the Corollaries 1 and
2 to our Lemma~\ref{CK} and Corollary~\ref{center}. Khovanov's Proposition 3
with Corollary 3 should be replaced by Lemma~\ref{capcup} and
Theorem~\ref{Endcapcup}. Corollary 4 corresponds to our Lemma~\ref{CK} again.
By repeating the arguments from \cite{Khovanovcob} the theorem follows.
\end{proof}

\begin{remark}
{\rm Although we have the invariant of the non-oriented tangles only up
  to shifts, the corresponding natural transformations would nevertheless
  satisfy the relations from~\cite{CRS}, since the only movie moves including
  $H$-moves are the moves $11$ to $14$ from Figure~\ref{fig:moves1114} and one
  can easily check that the corresponding natural transformations agree up to
  scalars, i.e. are well-defined.}
\end{remark}

\subsection{A conjectural connection with Khovanov homology}
In \cite{Khovanov} and \cite{Khotangles}, M. Khovanov introduced a homological
tangle and link invariant, now known as Khovanov homology. To any
$(2n,2m)$-tangle diagram he associated a certain complex of graded
$(\cH_n,\cH_m)$-bimodules for some combinatorially defined algebras
$\cH_n$, $\cH_m$. He proved that taking the graded cohomology groups defines an invariant
of tangles and links. As already mention in \cite{Braden}, there is a
connection between the algebras $A^\p$  from Section~\ref{Section1} where $\mg=\mathfrak{sl}_{2n}$ and $\p$
is the parabolic subalgebra given by all $(n,n)$-upper block matrices and
Khovanov's algebra $\cH_n$.

On the other hand, to each $(2n,2m)$-tangle diagram $t$ we associated a functor
$F=\Phi^{or}(t):\fD^b_{per}\big(\gMOF\text{-}A_{2n}\big)\longrightarrow\fD^b_{per}\big(\gMOF\text{-}A_{2m}\big)$
which can be described by tensoring with a complex $X_F$ of
graded $(A_{2m},A_{2n})$-bimodules. Let $F'$ denote the restriction of $F$ to a
functor
$F':\fD^b_{per}\big(\gMOF\text{-}A_{2n}^n\big)\longrightarrow\fD^b_{per}\big(\gMOF\text{-}A_{2m}^m\big)$
and let $X_F'$ be the associated complex of bimodules.
From Theorem~\ref{Tcomb} we know that we do not loose any information if we
restrict the functors to the categories of projective-tilting modules. Let
$\check{X}_F$ denote the corresponding complex of graded
$(D_{2m}^m,D_{2n}^n)$-bimodules. Let ${\bf H}^\bullet(F)$ denote the graded cohomology
of $\check{X}_F$.

The following conjecture relates the functorial invariant with Khovanov's
invariant:

\begin{conjecture}
\label{Khovlink}
  The homological tangle invariant $t\mapsto {\bf H}^\bullet(\check{X}_{\Phi^{or}(t)})$ is
  Khovanov's invariant.
\end{conjecture}

We would like to illustrate this on an

\begin{ex}
\label{itfits}
{\rm
  Let us consider the $2$-category of oriented $(2,2)$-tangles. It has a
  single object, $2$, and the functorial invariant $\Phi^{or}$ from
  Theorem~\ref{2catorient} assigns to each oriented $(2,2)$ tangle $t$ an
  endofunctor of the bounded derived category of $\MOF$-$A_2$ and, via
  restriction,  of
  $\MOF$-$A_2^1$. Recall that
  $\MOF$-$A_2^1\cong\cO_0(\mathfrak{sl}_2)$. From Example~\ref{ex2} we have
  $D_2^2\cong\mC[x]/(x^2)=\cA\langle 1\rangle$. Set $C=D_2^2$. Let us describe $\check{X}_{\Phi^{or}(t)}$ explicitly in terms of
  graded $C$-bimodules. To the flat diagram with two vertical strands we associated the identity functor which
  is given by tensoring with the graded $(C,C)$-bimodule $C$. To the other
  flat (cap-cup) diagram we associated the translation functor
  $\tilde{\theta}_1$. Via the functor $\mV^\mb$ from Corollary~\ref{faithful}
  this becomes tensoring with $C\otimes_\mC C\langle -1\rangle$ (see
  \cite{Sperv}). To the left twisted curls depicted in
  Figure~\ref{fig:oriented1} we associated the functor given by tensoring with the complex of graded $C$-bimodules
  \begin{eqnarray*}
\bigg(\cdots \rightarrow 0\rightarrow C\otimes_\mC C\langle -1\rangle
\stackrel{m}{\longrightarrow}  C\langle -1)\rightarrow 0\cdots\bigg)\langle
    -1\rangle\langle 3\rangle [1],
\end{eqnarray*}
which is
\begin{eqnarray}
\label{left}
\cdots\rightarrow 0\rightarrow C\otimes_\mC C\langle
1\rangle\stackrel{m}{\longrightarrow} C\langle 1\rangle\rightarrow 0\rightarrow\cdots,
  \end{eqnarray}
where $m$ is the multiplication map and $C\langle1\rangle$ is concentrated in
position zero of the complex (see also \cite[Section 8]{StDuke}). To the right twisted curl depicted in Figure~\ref{fig:oriented1} we associated the functor given by tensoring with the complex
  \begin{eqnarray*}
\bigg(\cdots \rightarrow 0\rightarrow  C\langle 1\rangle\stackrel{\Delta}\longrightarrow  C\otimes_\mC
C\langle -1\rangle \rightarrow0\cdots\bigg)\langle
    1\rangle\langle -3\rangle [-1],
\end{eqnarray*}
where $\Delta$ is mapping $1$ to $X\otimes 1+X\otimes 1$. Hence we get the complex
\begin{eqnarray}
\label{right}
\cdots \rightarrow 0\rightarrow  C\langle -1\rangle\stackrel{\Delta}\longrightarrow  C\otimes_\mC
C\langle -3\rangle \rightarrow0\cdots,
  \end{eqnarray}
where $C\langle -1\rangle$ is concentrated in position zero of the complex. 

On the other hand Khovanov associated to the flat diagram with two vertical
strands the algebra $\cA\langle1\rangle=C$ and to the other flat tangle the
$C$-bimodule $\cA\otimes\cA\langle 1\rangle\cong C\otimes_\mC
C\langle-1\rangle$. To the left twisted curls depicted in
Figure~\ref{fig:oriented1} corresponds the complex
\small
\begin{eqnarray*}
  \begin{array}[ht]{ccccccc}
&\bigg(\cdots\rightarrow 0\rightarrow &\cA\otimes\cA\langle
1\rangle&\stackrel{m}{\longrightarrow} &\cA\langle 1-1\rangle\rightarrow&
0\rightarrow\cdots\bigg)\langle 2\rangle [-1],\\
=\bigg(\cdots\rightarrow 0\rightarrow&\cA\otimes\cA\langle
3\rangle\stackrel{m}{\longrightarrow}& \cA\langle 2\rangle&\rightarrow
&0\rightarrow\cdots\bigg)
\end{array}
  \end{eqnarray*}
\normalsize
where $\cA\langle 2\rangle\cong C\langle 1\rangle$ is concentrated in position
zero.
Hence the complex coincides with \eqref{left}. The complex associated in
\cite{Khotangles} to the right twisted curl depicted in
Figure~\ref{fig:oriented2} is exactly the complex \eqref{right}.}
\end{ex}
Therefore, in the example of $(2,2)$-tangles, Khovanov homology is nothing else than the homology of the
functorial invariant from \ref{2catorient} restricted to the underlying category
given by projective-tilting modules.

\section{Kac-Moody algebras}
In this section we consider projective (respectively tilting) functors for
symmetrizable Kac-Moody algebras. Given any regular block outside the critical
hyperplanes, the main result is
the classification theorem of these projective endofunctors - which can be viewed as a categorification of the group algebra of the integral
Weyl group associated to the block. The result was
conjectured in \cite{FMengl}, where the authors defined tilting functors via
the so-called Kazhdan-Lusztig tensoring. Since this construction is highly non-trivial we
decided to work with an alternative definition instead. Our approach is based on the
translation functors introduced and described in \cite{Neidhardt1}, \cite{Neidhardt2} and \cite{CT}.\\

Let now $\mg$ be a symmetrizable Kac-Moody algebra (in the sense of
\cite{Kacbook}) over the complex numbers with a fixed triangular
decomposition $\mg=\mn^{-}\oplus\mh\oplus\mn=\mn^-\oplus\mb$. Let
$\cU=\cU(\mg)$, $B=\cU(\mb)$ and $S=\cU(\mh)$ denote the corresponding
universal enveloping algebras. Let $\cW$ be the Weyl group. Let $T$ be a {\it deformation algebra},
i.e an associative commutative noetherian unital $S$-algebra. Important
examples will be $T=\mC\cong S/(\mh)$ and $T=S_{(0)}$ which is the localisation
of $S$ at $(\mh)$.  Let $\op{t}:S\rightarrow T$ be the
 morphism defining the $S$-structure on $T$.\\
Since $\mg$ is symmetrizable we have a non-degenerated symmetric
bilinear form on $\mg$ inducing a non-degenerated bilinear form
$(\bullet,\bullet):\mh^\ast\times\mh^\ast\rightarrow\mC$ on $\mh$
via restriction and dualising. This extends $T$-linearly to a
bilinear form $(\bullet, \bullet)_T:\mh^\ast\otimes T\times
\mh^\ast\otimes T\rightarrow T$. We consider the restriction $\op{t}_\mh$ of
$\op{t}$ to $\mh$ as an element of $\HOM_\mC(\mh,T)=\mh^\ast\otimes T$ and
define $h_\mu=(\mu,\op{t}_\mh)_T$  for any $\mu\in\mh^*\subset\mh^*\otimes T$. Let
$(\bullet,\bullet)_\mC$ be the specialisation of $(\bullet,\bullet)_T$. A
weight $\la\in\mh^\ast$ lies {\it outside the critical hyperplanes} if
$(\la+\rho+\op{t}_\mh,\beta)_\mC\notin\mZ(\beta,\beta)_\mC$ for any imaginary
root $\beta$. We fix $\rho\in\mh^\ast$ such that
$(\rho,\alpha)=1$ for any simple root $\alpha$. The dot-action of
the Weyl group on $\mh^\ast$ is defined as in the finite
dimensional case. In the following we assume the reader to be
familiar with the definition and results of \cite{CT} and
\cite{Fiebigcomb}.

\subsection{The deformed category $\cO_T$}
\label{deformed}
We consider the {\it $T$-deformed category $\cO$} denoted by $\cO_T$. It is
the full subcategory of the category of $\cU\otimes T$-modules given by locally $B\otimes T$-finite objects having a
weight decomposition $M=\bigoplus_{\la\in\mh^*}M_\la$, such that
\begin{eqnarray*}
  M_\la=\{m\in M\mid hm=(\la+\op{t})(h)m, \text{ for all }h\in\mh\},
\end{eqnarray*}
where the $(\la+\op{t})(h)$ are considered as elements of $T$, see
   e.g. \cite[Section 2.1]{CT}. Set
   $\mg_T=\mg\otimes T$. 
For $\la\in\mh^\ast$ the composition $B\rightarrow S\stackrel{\la+\op{t}}\lra
   T$ defines a $B$-structure on $T$ which commutes with the usual left $T$-action. The resulting $B\otimes
   T$-structure on $T$ will be denoted by $T_\la$. We denote by
   $\Delta_T(\la)\in\cO_T$ the {\it $T$-deformed Verma
     module} $\Delta_T(\la)=\cU\otimes_B T_\la$ with highest weight
 $\la$.  If $T$ is local with maximal ideal $\mm$ or $T=\mC$ then the isomorphism classes
   $L_T(\la)$ of
   simple objects in $\cO_T$ are parametrised by elements $\la\in\mh^\ast$ by
   taking their highest weights. Under the specialisation
   functor $T/\mm\otimes_T\bullet$ the simple objects in $\cO_T$ become the
   simple objects in the ordinary BGG-category $\cO_{T/\mm}$ (where the modules need not to be finitely generated). For details we refer to \cite[Proposition 2.1]{CT}.

\subsection{Blocks outside the critical hyperplanes}
\label{transl}
Let now $T$ be a {\it local} deformation algebra with maximal ideal $\mm$ and
residue field isomorphic to $\mC$. There is a block decomposition (\cite[Section 2.4]{CT})
\begin{eqnarray*}
\cO_T=\Pi_\Lambda\cO_{T,\Lambda}
\end{eqnarray*}
indexed by certain
sets $\Lambda$  of weights (the highest weights of the simple objects in the
block). The block $\cO_{T,\Lambda}$ is the full subcategory of $O_{T}$ given by
all objects where the highest weights of any subquotient are contained in
$\Lambda$. If all the weights of
$\Lambda$ are outside the critical hyperplanes then we call $\cO_{T,\Lambda}$
{\it a block outside the critical hyperplanes}. If moreover $T$ is a local
deformation domain then $\Lambda$ is the
dot-orbit of any of its elements under the corresponding {\it integral} Weyl
group $\cW_T(\Lambda)$ (which is generated by reflections corresponding to
real roots only, \cite[page 699]{CT}). If $T\rightarrow T'$ is a homomorphism of local deformation
algebra domains, then the base change functor $T'\otimes_T\bullet $ maps the
blocks of $\cO_T$ to the blocks of $\cO_{T'}$ (\cite[Corollary 2.10]{CT}). On
the other hand, for $T'=T/\mm$ the resulting decomposition is
exactly the block decomposition from \cite[Theorem 2]{KacKaz}. In case $T\rightarrow T'$
is a homomorphism of deformation algebras we denote by  $\cO_{T',\Lambda}$ the
image  of a block $\cO_{T,\Lambda}$ under the base change functor. In
general, this is not a block, but a direct sum of blocks (\cite[Lemma 2.9,
Corollary 2.10]{CT}). For any $\cO_{T,\Lambda}$ we denote by $\cM_{T,\Lambda}$ the full subcategory
of $\cO_{T,\Lambda}$ given by all modules having a finite Verma flag.

\subsection{Translation through walls}
Let $T$ be any local deformation algebra domain. Let $\cO_{T,\Lambda}$,
$\cO_{T,\Lambda'}$ be two blocks, outside the critical hyperplanes. Let $\la\in\Lambda$, $\la'\in\Lambda'$. Assume
that
\begin{enumerate}[(1)]
\item $\la-\la'$ is integral and there is a dominant weight $\nu$ in the
$\cW$-orbit of $\la-\la'$ (in particular $\cW_T(\Lambda)=\cW_T(\Lambda')$),
\item $\la$ and $\la'$ lie in the closure of the same Weyl chamber,
\item under the dot-action, the $\cW_T(\Lambda)$-stabiliser of $\la$ is
  contained in the $\cW_T(\Lambda')$-stabiliser of $\la'$,
\end{enumerate}
then there are the translation functors {\it onto the wall} and {\it out of the
  wall}
\begin{eqnarray*}
  \theta_{on}: \cM_{T,\Lambda}\rightarrow\cM_{T,\Lambda'}&&
  \theta_{out}: \cM_{T,\Lambda'}\rightarrow\cM_{T,\Lambda}
\end{eqnarray*}
as defined in \cite{CT}. The definition of $\theta_{on}$ is completely
analogous to the finite dimensional case, and given by tensoring with the
simple integrable highest weight module $L(\nu)$ and
projection onto the block. The functor $\theta_{out}$ is defined by taking a
limit of certain truncations of the functor tensoring with the restricted dual
of $L(\nu)$ (see \cite[4.1]{CT}).
If $\la'$ lies exactly on the $s$-wall, i.e. the $\cW_T(\Lambda')$-stabiliser
of $\la'$ is $\{e,s\}$ and $\la$ is regular, then we have the translation
functor $\theta_s=\theta_{out}\theta_{on}$ through the $s$-wall. Note that
translation functors commute with base change in the sense of
\cite[5.2]{CT}. In particular, it gives translation functors for any
$S_{(0)}$-algebra $T'$.

\subsection{The fake antidominant projective module}
Let $\cO_{T,\Lambda}\subset\cO_T$ denote a regular block where all the
  corresponding weights are outside the critical hyperplanes. In general,
  $\cO_{T,\Lambda}$ does not have enough projectives. In particular, in the Kac-Moody case, Soergel's famous
  ``antidominant projective module'' does not need to exist. However,
  there are (\cite[Theorem 2.7]{CT}) enough projectives in the truncated categories
  $\cO_{T,\Lambda}^{\nu}$ given by all objects in $\cO_{T,\Lambda}$ whose
  weights are $\leq\nu$ (with $\nu\in\mh^\ast$ fixed). If $\Lambda$ contains
  an antidominant weight $\lambda$, we denote by
$P_T^n(\la)$ the projective cover of the simple module $L_T(\la)$ in
$\cO_{T,\Lambda}^{\la+n\chi}$, where $\chi$ is the sum of all simple
roots. We choose a compatible system of surjections $p_{m,n}:P^{m}_T(\la)\surj
  P^n_T(\la)$ for $m\geq n$ and denote $P^\infty_T(\la)=\varprojlim P^n(\la)$. By definition, there is a canonical isomorphism
$$\HOM_{\mg_T}(P^\infty_T(\la),\varprojlim
  P^n_T(\la))\cong\varprojlim\HOM_\mg(P^\infty_T(\la)P^n_T(\la)).$$ Since the largest quotient of $P_T^\infty(\la)$ contained
  in $\cO^{\la+n\chi}$ is $P^n_T(\la)$, the induced inclusion
  \begin{eqnarray*}
    \HOM_{\mg_T}(P^n_T(\la),P^n_T(\la))\inj\HOM_{\mg_T}(P^\infty_T(\la),P^n_T(\la))
  \end{eqnarray*}
is in fact an isomorphism. Therefore, we have canonically
\begin{eqnarray}
\label{einschub}
\varprojlim\END_{\mg_T}(P_T^n(\la))&\cong&\varprojlim\HOM_{\mg_T}(P^\infty_T(\la),P^n_T(\la))\nonumber\\
&\cong&\HOM_{\mg_T}(P^\infty_T,
  \varprojlim P^n_T(\la))\nonumber\\
&=&\END_{\mg_T}(P_T^\infty(\la)).
\end{eqnarray}
If $F:\MTL\rightarrow\MTL$ is a
composition of translations through walls we set $F\Pinf=\varprojlim F
P^n_T(\la)$, where the defining maps are the $F(p_{m,n})$ for $m\geq n$.

\subsection{Blocks containing an antidominant weight}
\label{general} The definition of translation functors through
walls is based on the existence of some dominant weight $\nu$ as
stated above. Such a dominant weight need not to exist in general
for any block $\cO_{T,\Lambda}$ outside the critical hyperplanes.
However, in these cases there will be an antidominant weight
$\la\in\Lambda$. We will study such blocks now. Let $T$ be a local deformation algebra and let
  $\cO_{T,\Lambda}$ be a regular block outside the critical hyperplanes. Let
  us assume $\Lambda$ to have an antidominant weight $\lambda$. Let
  $\tau(\Lambda)=\{-2\rho-\lambda\mid\la\in\Lambda\}$. In particular,
  $\tau(\Lambda)$ contains a dominant weight, and therefore $\cO_{\Lambda,T}$
  has enough projectives. We consider the (Chevalley-)anti-automorphism $\sigma$
  interchanging root spaces $\mg_\alpha$ and
  $\mg_{-\alpha}$ and the principal anti-automorphism $\gamma:S\rightarrow S$,
  $h\mapsto -h$ for any $h\in\mh$. For a $\mg$-module $M$, we denote by
  $_{}^\sigma M$ the space $M$ but with $\sigma$-twisted $\mg$-action. Likewise,
  given an $S$-module $N$, the symbol $_{}^{\gamma}N$ denotes $N$ with
  $\gamma$-twisted action of $S$. The {\it semi-infinite bimodule} $S_{-2\rho}$ with respect to the
  semi-infinite character $-2\rho$ is defined in \cite{SoKac} (and relies on
  the work of Arkhipov, Frenkel and Voronov). It provides an
  equivalence of categories
  \begin{eqnarray}
  \label{tiltingequiv}
    \tau:\cM_{T,\Lambda}\iso\cM_{_{}^\gamma T,\tau(\Lambda)}^{\op{opp}}
  \end{eqnarray}
which maps short exact sequences to such (\cite{SoKac}). More precisely, the equivalence is constructed as follows: We consider the
functor
\begin{eqnarray*}
  \cT_{-2\rho}: \quad M\mapsto \quad^\sigma\big((S_{-2\rho}\otimes_\cU M)^\oast\big)
\end{eqnarray*}
from the category of $\mg_T$-modules whose weight spaces
are free $T$-modules of finite rank. Hence $(S_{-2\rho}\otimes_\cU
M)^\oast:=\bigoplus_{\la\in\mh\ast}\HOM_T((S_{-2\rho}\otimes_\cU
M)_\la,T)$ is again a free $T$-module. The natural right
$\mg$-action becomes a left $\mg$-action after twisting with
$\sigma$. Hence $\cT_{-2\rho}(M)$ is a $\mg_T$-module. The
arguments in \cite{SoKac} show that $\cT_{-2\rho}$ restricts to an
equivalence \eqref{tiltingequiv}, called the {\it tilting equivalence},
sending $\triangle_T(\mu)$ to $\triangle_{^\gamma T}(\tau(\mu))$
for any $\mu\in\Lambda$. Note that if $T=S_{(0)}$ or $T=\mC$, then
$_{}^\gamma T\cong T$ as $S$-modules. Therefore, if
$P\in\cO_{T,\tau(\Lambda)}$ is projective and finitely generated,
hence in $\cM_{T,\Lambda}$,  then
$\tau(P)\in\cM_{T,\Lambda}^{\op{opp}}$ is {\it tilting}, i.e. it
has a finite Verma flag and
$$\EXT^1_{\cO_{T,\Lambda}}(\Delta(\mu),\tau(P))=\EXT^1_{\cO_{{^\gamma}T},\tau(\Lambda)}(P,\Delta(\tau(\mu))=0$$
for any $\mu\in\Lambda$. It follows directly from the definitions
that $\cT_{-2\rho}$ commutes with base change. Via these tilting
equivalences the translation functors $\theta_s$ from
Section~\ref{transl} give rise to translation functors $\theta_s$
for blocks containing an antidominant weight.

In the following an {\it antidominant block} means a regular block $\cO_{T,\Lambda}$ outside the critical hyperplanes such that $\Lambda$ contains an
antidominant weight which we denote by $\la$.

We first need an analogue of Proposition~\ref{Irving}, namely that the simple objects
occurring in the socle of a Verma module in $\cO^\p$ are of maximal
Gelfand-Kirillov.

\begin{lemma}
\label{A}
  Let $\cO_{\mC,\Lambda}$ be an antidominant block. Let $X\in\cM_{\mC,\Lambda}$. Then
  $\HOM_{\mg}\big(P^\infty_\mC(\la),X'\big)\not=0$ for any submodule $X'$ of $X$.
\end{lemma}
\begin{proof}
  As a submodule of a module with Verma flag, $X'$ contains a Verma module,
  hence also $\Delta_\mC(\la)$ by \cite[Theorem 3.10]{Fiebigcomb}.
\end{proof}

The following result is again well-known in the finite dimensional situation
(see e.g. \cite[4.13 (1)]{Ja2}):
\begin{lemma}
\label{B}
   Let $T$ be a local deformation algebra domain or $T=\mC$. Let
  $\cO_{T,\Lambda}$ be an antidominant block. Then $P^n_T(\la)$ has a finite Verma flag and for any
  $\mu\in\Lambda$ we have $(P_T^n(\la)\::\:\Delta_T(\mu))=1$ if $n\gg0$.
\end{lemma}
\begin{proof}
  The existence of a finite Verma flag is given by \cite[Theorem 2.7]{CT}. The
  multiplicity formula is \cite[Lemma 3.8]{Fiebigcomb}.
\end{proof}

\subsection{The centre}
Recall from Section~\ref{Section1} that the centre $C$ of the
algebra $A^\mb$, or equivalently the centre of the category
$\cO_0^\mb$, is isomorphic to the endomorphism ring of the full
tilting module $T$, which is Soergel's "antidominant projective
module", and the classification theorem of projective functors can
be obtained by considering certain special $C$-bimodules. We would
like to generalise this approach to the Kac-Moody case. Therefore
we first describe the centre of the deformed categories
$\cO_{T,\Lambda}$ (generalising the main result of \cite{CT}).

\begin{lemma}
\label{thetaP}
  Let $T=S_{(0)}$ or $T_\mC$. Let $\cO_{T,\Lambda}$ be an antidominant
  block and $s\in\cW(\Lambda)$ a simple reflection. Then $\theta_s \Pinf\cong\Pinf\oplus\Pinf$.
\end{lemma}

\begin{proof}
  By definition of the translation functors, $\theta_s P^n_T(\la)$ is
  projective in a suitable truncation $\cO_{T,\Lambda}'$ of
  $\cO_{T,\Lambda}$. Since $T=S_{(0)}$ or $T=\mC$ the functor $\theta_s $ is
  self-adjoint (\cite[Corollary 5.10]{CT}). Since $\dim_\mC\HOM_{\mg_T}(P^n_\mC(\la), \theta_s\Delta_\mC(\la))=2$ by
  \cite[Corollary 5.10, Proposition 4.1]{CT} and \cite[Remark 3.9
  (2)]{Fiebigcomb} we get that the indecomposable
  cover of $\Delta_\mC(\la)\in \cO_{\mC,\Lambda}'$ occurs with multiplicity
  $2$ as a direct summand of $\theta_s P^n_\mC(\la)$.  Therefore, in
  $\theta_s P^n_{S_{(0)}}(\la)$, the indecomposable
  cover of $\Delta_{S_{(0)}}(\la)$ occurs with multiplicity
  $2$ as well (\cite[Proposition 2.6, Lemma 5.4]{CT}), since the blocks in
  $\cO_\mC$ are the specialisations of the blocks of $\cO_{S_{(0)}}$.
If another direct summand occurs, than it has to occur for any
  $m>n$. However, Lemma~\ref{B} says that any Verma module occurs in
  $\theta_s P_T^\infty(\la)$ once, and \cite[Proposition 4.1]{CT}
  implies that any Verma module occurs in $\theta_s\Pinf$ with multiplicity
  $2$. It follows $\theta_s\Pinf\cong \Pinf\oplus\Pinf$.
\end{proof}

\begin{remark}
\label{baseanti}
  {\rm If $T$ is any local deformation algebra then $\theta_s
  \Pinf$ is isomorphic to a finite direct sum of fake antidominant projective
  modules $P_T^\infty(\la')$, each occurring with finite multiplicity. (\cite[Lemma
  5.4]{CT}, \cite[Remark 3.9]{Fiebigcomb}).}
\end{remark}

\begin{theorem}
\label{centerKac}
Let  $\cO_{S_{(0)},\Lambda}$ be an antidominant block. Let $T=\mC$ or
${S_{(0)}}\rightarrow T$ be a morphism of local deformation algebras. There is a natural isomorphism
  \begin{eqnarray*}
    \cI=\cI_{T,\Lambda}:\quad\END(\ID_{\cO_{T,\Lambda}})&\cong&\END_{\mg_T}(T\otimes_{S(0)}
    P^\infty_{S(0)}(\la))
  \end{eqnarray*}
\end{theorem}
\begin{proof}
Using the identifications \eqref{einschub} we claim that
$\phi\mapsto\phi_{\infty}:=\{\phi_{P^n_T(\la)}\}$ defines the required
isomorphism.

  {\it Injectivity:}
  If $T=\mC$ then the injectivity follows directly from Lemma~\ref{A}. For
  the general case we have to work more. Assume $\phi_{\infty}=0$. Assume now $T={S_{(0)}}$ and let
  $X\in\OTL$ be a tilting module (having a finite Verma flag). Let us assume
  $X$ is indecomposable. By
  \cite[Section 4]{Fiebigcomb} $X$ is a direct summand of some
  $F\Delta_T(\la)$, where $F$ is a composition of translations through walls.
  Since $P^n_T(\la)$ has a finite Verma flag with $\Delta_T(\la)$ at the top,
  and $F$ preserves short exact sequences of modules with finite Verma flag (\cite[Proposition 4.1]{CT}),
  $FP^n_T(\la)$ surjects onto $X$. By Lemma~\ref{thetaP}, a (finite) direct
  sum  of $\Pinf$ surjects onto $X$. Hence $\phi_X=0$. Since $\OTL$ is
generated by modules having a finite Verma flag (\cite[Theorem 2.7]{CT}),
$\cI$ is injective. The general case follows by base change using Remark~\ref{baseanti}.\\
{\it Surjectivity:}
Via restriction the modules
$\{\cZ^n_{T,\Lambda}:=\END(\ID_{\cO_{T,\Lambda}^{\la+n\chi}})\}_{n\in\mN}$ form
a projective system and
$\cZ_{T,\Lambda}:=\END(\ID_{\cO_{T,\Lambda}})\cong\varprojlim\cZ^n_{T,\Lambda}$.
It is sufficient to show that
$\END(\ID_{\cO_{T,\Lambda}^{\la+n\chi}})\inj\END_{\mg_T}(T\otimes_{S(0)}P^n_{S(0)}(\la))$,
  $\phi\mapsto \phi_{P^n_T(\la)}$ defines an isomorphism for any $n\in\mN$.
If $T=S_{(0)}$ then the statement is true by \cite[Lemma
  3.12]{Fiebigcomb}. Let now $Q=\op{Quot}(S_{(0)})$ be the quotient field of $S_{(0)}$. Since there are no
extensions between Verma modules, there are canonical isomorphisms
\begin{eqnarray*}
  \END\big(\ID_{\cO^{\la+n\chi}_{Q,\Lambda}}\big)\cong\prod_{\la+n\chi\geq\mu\in\Lambda}\END_{\mg_Q}\big(\Delta_Q(\mu)\big)\cong\prod_{\la+n\chi\geq\mu\in\Lambda}Q.
\end{eqnarray*}
On the other hand
\begin{eqnarray*}
\END_{\mg_Q}\big(Q\otimes_{S(0)}P_{S(0)}^n(\la)\big)&\cong&
\END_{\mg_Q}\Big(\prod_{\la+n\chi\geq\mu\in\Lambda}\Delta_Q(\mu)\Big)\\
&\cong&\prod_{\la+n\chi\geq\mu\in\Lambda}
\END_{\mg_Q}\big(\Delta_Q(\mu)\big)
\end{eqnarray*}
by Lemma~\ref{B}. Since the maps in question are morphisms of free $Q$-modules
  of finite rank, $\cI$ becomes an isomorphism $\cI_{Q,\Lambda}$ at the
  generic point. We have to show that it also
  becomes an isomorphism $\cI_{\mC,\Lambda}$ at the closed point to get the
  surjectivity in general. For this let
  $f\in\END_{\mg_\mC}(\mC\otimes_{S(0)}P^\infty_{S(0)})$. 
  Recall that specialisation defines an isomorphism
  $\END_{\mg_\mC}(\mC\otimes_{S(0)}P^\infty_{S(0)}(\la))
\cong\mC\otimes_{S(0)}\END_{\mg_{S_{(0)}}}(P^\infty_{S(0)}(\la))$
  (\cite[Proposition 2.4]{CT} or \cite[Remark 3.9]{Fiebigcomb})
  under which $f$ corresponds to some $1\otimes g$. Since
  $\cI_{S_{(0)},\Lambda}$ is an isomorphism, there exists some $\phi\in
  \END(\ID_{\cO_{S_{(0)},\Lambda}})$ such that
  $\cI_{S_{(0)},\Lambda}(\phi)=g$. Then $1\otimes \phi\in
  \END(\ID_{\cO_\mC,\Lambda})$ (\cite[Proposition 3.1]{CT}) and by definition
of $\cI$ we get  $\cI_{\mC,\Lambda}(1\otimes\phi)=\cI_{S_{(0)},\Lambda}(1\otimes \phi)=1\otimes g=f$. This means
$\cI_{\mC,\Lambda}$ is surjective. The statement of the theorem follows.
\end{proof}

We denote $\cZ_{T,\Lambda}=\END(\ID_{\OTL})$, the centre of $\cO_{T,\Lambda}$.
\begin{cor}
\label{corcenter}
Let $\cO_{S_{(0)},\Lambda}$ be an antidominant block. Let $T=\mC$ or
 $S_{(0)}\rightarrow T$ be a morphism of local deformation algebras. Then
 there is a natural isomorphism $$\cZ_{T,\Lambda}\cong T\otimes_{S(0)}\cZ_{S_{(0)},\Lambda}.$$
\end{cor}
\begin{proof}
  Let $n\in\mN$. There are isomorphisms of rings
\begin{eqnarray*}
 T\otimes_{S(0)}\END_{\mg_{S(0)}}\big(P^n_{S(0)}(\la)\big)
\cong\END_{\mg_T}\big( T\otimes_{S(0)}P^n_{S(0)}(\la)\big)
\end{eqnarray*}
by \cite[Proposition 2.4]{CT} and \cite[Proposition 2.6]{CT}.
Taking limits the theorem gives the statement.
\end{proof}

\subsection{The structure theorem}
Let $\cO_{S_{(0)},\Lambda}$ be an antidominant block. Let $S_{(0)}\rightarrow
T$ be a morphism of local deformation algebras. We consider the functor $\mV_T=\HOM_{\mg_T}(T\otimes_{S(0)}\PRinf,\bullet):\OTL\rightarrow
  \MOD\text{-}\cZ_{T,\Lambda}$. Note that with the assumptions of
  Corollary~\ref{corcenter} and the formulas~\eqref{einschub} there are canonical isomorphisms of
$\cZ_{T,\Lambda}$-modules
\begin{eqnarray}
\label{VLIMS}
\begin{array}[tbc]{llllll}
\mV_{S(0)}\varprojlim P_{S(0)}^n(\la)&\cong&\varprojlim\mV_{S(0)} P_{S(0)}^n
&\cong&
\varprojlim\END_{\mg_{S(0)}}(P_{S(0)}^n(\la))\\
&\cong&\END_{\mg_{S(0)}}(P_{S(0)}^\infty)&\cong&\cZ_{{S(0)},\Lambda}.
\end{array}
\end{eqnarray}
Using Corollary~\ref{corcenter}, specialisation gives an isomorphism 
\begin{eqnarray}
\begin{array}[tbc]{llllll}
\mV_T\varprojlim (T\otimes_{S(0)}P_{S(0)}^n(\la))
&\cong&
\varprojlim\mV_{S(0)} (T\otimes_{S(0)} P_{S(0)}^n)\\
&\cong&
T\otimes_{S(0)}(\varprojlim\END_{\mg_{S(0)}}(P_{S(0)}^n(\la)))\\
&\cong&
T\otimes_{S(0)}\END_{\mg_{S(0)}}(P_{S(0)}^\infty)\\
&\cong&T\otimes_{S(0)}\cZ_{{S(0)},\Lambda}\\
&\cong&
\cZ_{T,\Lambda}.
\end{array}  
\end{eqnarray}

Recall that the deformed Verma modules in an antidominant block
$\cO_{S_{(0)},\Lambda}$ are exactly the ones with highest weight in
$\Lambda=\cW(\Lambda)\cdot\la$. In \cite[Theorem 3.6]{CT} one can find an
explicit description of the centre 
$\cZ_{S_{(0)},\Lambda}$ as a subring of $\Pi_{w\in\cW(\Lambda)}S_{(0)}$ 
by looking at the acting on each deformed Verma module. Together with
Corollary~\ref{corcenter} we get a concrete description of
$\cZ_{\mC,\Lambda}$. From this description we also get a natural right action
of $\cW(\Lambda)$ on $\cZ:=\cZ_{S_{(0)},\Lambda}$. For any simple reflection
$s\in\cW(\Lambda)$ we denote by $\cZ^s$ its invariants and get the following
important result

\begin{lemma}
\label{freeness}
$\cZ$ is a free $\cZ^s$-module and gives rise to a self-adjoint functor
functor $$\bullet_{\cZ^s}\cZ:\MOF\text{-}\cZ\rightarrow\MOF\text{-}\cZ.$$ Moreover,
there is an isomorphism of functors $\mV_{S(0)}\theta_s(\bullet)\cong \mV_{S(0)}(\bullet)\otimes_{\cZ^s}\cZ$.
\end{lemma}

\begin{proof}
Set $R=S_{(0)}$. The following argument is due to
P. Fiebig: Let $\alpha$ be the simple root corresponding to $s$. For $w\in\cW$
such that $w\cdot\la\in\Lambda$ we set $z_w=h_{w(\alpha)}$. In particular,
$z_{w}=-z_{ws}$ and also $\{z_w\}\in\cZ$ by \cite[Theorem 3.6]{CT}. Let now
$a=\{a_w\}$ be an element of $\cZ$. We claim that there are
uniquely defined elements $c_+, c_-\in\cZ^s$ such that $a=c_++c_-z$ and (hence
$1$, $z$ the desired basis). 
From the general equality $ws=s_{w{\alpha}}w$, we get
$a_{ws}=a_{s_{w(\alpha)}}\equiv a_w\mod z_w$ by \cite[Theorem
  3.6]{CT}. Let now 
  $r_1$, $r_2\in R$ such that
  $r_1\equiv r_2 \mod h$ for some $0\not=h\in\mh$.  Then the equation
  $(r_1,r_2)=c_+(1,1)+c_-(h,-h)$ has the unique solution
  $c_\mp=\frac{1}{2h}(r_1\mp r_2)$. By \cite[Theorem 3.6]{CT} the elements
  $x_+$, $x_-$ exist and are contained in $\cZ$, because $a\in\cZ^s$.  The
  self-adjointness follows then as in \cite[Proposition 5.10]{Sbimods}. The
  last statement is \cite[Theorem 4.1]{Fiebigcomb}.
\end{proof}

\begin{lemma}
\label{clear} Let $\cO_{S_{(0)},\Lambda}$ be an antidominant
block. Let $T=\mC$ or let $S_{(0)}\rightarrow T$ be a morphism of local
deformation algebras. Then $\mV_T$ induces an isomorphism
  \begin{eqnarray*}
    \Phi: \HOM_{\mg_T}(X,Y)
\cong
\HOM_{\cZ_{T,\Lambda}}(\mV_T X,\mV_T Y)
  \end{eqnarray*}
for $X=T\otimes_{S(0)}\PRinf$ and any $Y\in\OTL$ or
$Y=T\otimes_{S(0)}F\PRinf$, where $F$ is a finite composition of translation functors through walls.
\end{lemma}
\begin{proof}
  Let first $Y\in\OTL$. By definition of $\mV_T$ we have $$\mV_T Y=\HOM_{\mg_T}(T\otimes _{S(0)}\PRinf,Y)\rightarrow \HOM_{\cZ_{T,\Lambda}}(\mV_T X,\mV_T
  Y)\cong\mV_TY,$$ the latter by evaluating $f$ at the identity morphism
  $\op{id}$. This composition is mapping $f\in\mV_T Y$ to
  $\mV_T(f)(\op{id})=f\circ\op{id}=f$, hence the middle arrow is an isomorphism
  and the statement follows for any $Y\in\OTL$. Then the lemma follows by
  taking limits.
\end{proof}

The following is a crucial refinement of \cite[Theorem 3.25]{Fiebigcomb}.
\begin{theorem}[Structure theorem]
\label{injectivity} Let $\cO_{S_{(0)},\Lambda}$ be an antidominant
block. Let $T=\mC$ or let
$S_{(0)}\rightarrow T$ be a morphism of local deformation algebras.
Let $X$, $Y\in\OTL$ be tilting modules. Then $\mV_T$ induces a
natural inclusion
  \begin{eqnarray}
\label{incl}
    \Phi: \HOM_{\mg_T}(X,Y)\inj\HOM_{\cZ_{T,\Lambda}}(\mV_T X,\mV_T Y)
  \end{eqnarray}
which is an isomorphism (at least) if $T=S{(0)}$ or $T=\mC$.
\end{theorem}

\begin{proof} Recall from \cite{Fiebigcomb} that any indecomposable tilting
module $X\in\cO_{T,\Lambda}$ is a direct summand of some $T\otimes _{S(0)}(F\Delta_{S(0)}(\la))$
for some composition $F$ of translations through walls. Hence it is enough to
show the statement for $X=T\otimes_{S(0)}(F\Delta_{S(0)}(\la))$.
Choose compatible surjections $P^n_{S(0)}(\la)\rightarrow\Delta_{S(0)}(\la)$ for
$n\geq 0$. Then $F P^n_{S(0)}(\la)\rightarrow F \Delta_{S(0)}(\la)$ is again
surjective (\cite[Proposition 4.1]{CT}) and so is the
specialisation  $T\otimes_{S(0)}(F P^n_{S(0)}(\la))\rightarrow T\otimes_{S(0)}
(F\Delta_{S(0)}(\la))$, and it stays surjective if we apply the
exact functor $\mV_T$. Thanks to Lemma~\ref{thetaP} it is enough to verify \eqref{incl} for $X=T\otimes_{S(0)}\PRinf$. This is
however Lemma~\ref{clear}. Let now $T=S_{(0)}$ or $T=\mC$. Because of the self-adjointness of translations through
walls \cite[Corollary 5.10, Proposition 3.11 (4)]{Fiebigcomb} and Lemma~\ref{freeness} we may assume
that $F$ is isomorphic to the identity, hence $X\cong\Delta_T(\la)$. Let $I$
be the maximal ideal of $\cZ_{T,\Lambda}=\END_{\mg_T}(\Pinf)$. 
Then $\HOM_{\mg_T}(X,Y)$ can be identified with the space $\{f\in
\HOM_{\mg_T}(\Pinf,Y)\mid f\circ g=0, g\in I\}$ and $\HOM_{\mg_T}(\mV_T X,\mV_T Y)$ can be identified with the space $\{f\in
\HOM_{\cZ_{T,\Lambda}}(\mV_T\Pinf,Y)\mid f\circ g=0, g\in I\}$. So, the 
theorem follows from Lemma~\ref{clear}.   
\end{proof}

\subsection{The combinatorics of natural transformations}
Our next step is to prove a generalisation of Theorem~\ref{Tcomb} for
Kac-Moody algebras. We
start with some preparations. From now on we fix an antidominant block
$\cO_{S_{(0)},\Lambda}$. Let $T=S(0)$ or $T=\mC$ and denote $P^m=P^m_T(\la)$ and
$P=P^\infty_T(\la)$. Note that, for any finite composition $F$ of translation
functors through walls, the module $FP$ has
naturally a $\cZ_{T,\Lambda}$-module structure, by definition of the centre. This action
commutes with the second left $\cZ_{T,\Lambda}$-structure given by $z.m=F(z)(m)$ for any
$z\in\cZ_{T,\Lambda}=\END_{\mg_T}(P)$, $m\in F P$. Via the
functor $\mV_T$ the first left action converts into the usual right $\cZ_{T,\Lambda}$-action
on $\mV_T F P$ and the second actions turns into the left $\cZ_{T,\Lambda}$-action
$z.f=F(z_P)\circ f$. We have a projective system
\begin{eqnarray*}
  \big\{H_m:=\HOM_{\mg_T}(P, FP^m (\la))=\mV_TFP^m(\la)\big\}_{n\in\mN}
\end{eqnarray*}
given by the maps $F(p_{m,n})\circ\cdot:H_m\rightarrow H_n$ for $m\geq n$. Set
$\mV_T FP (\la)=\HOM_{\mg_T}\big(P,\varprojlim
FP^n(\la)\big)\cong\varprojlim\HOM_{\mg_T}(P(\la),FP^n(\la))$. We get the following generalisation of Theorem~\ref{Tcomb}:
\begin{theorem}
\label{homprojfuncKM}
Let $T=S_{(0)}$ or $T=\mC$ and let $\cO_{T,\Lambda}$ be an antidominant block. 
Let $F$, $G:\MTL\rightarrow\MTL$ be finite compositions of translations through
  walls. There is a natural isomorphism of vector spaces (even of rings if $F=G$)
  \begin{eqnarray*}
    \HOM(F,G)\cong\HOM_{\cZ_{T,\Lambda}\text{-}\MOD\text{-}\cZ_{T,\Lambda}}(\mV_T
    FP^\infty(\la),\mV_T GP^\infty(\la))
  \end{eqnarray*}
\end{theorem}

\begin{proof} We claim that $\phi\mapsto a^\phi\in\fH$
where $(a^\phi)_n(\{g_j\}_{j\in\mN})=\phi_n\circ g_n$ with
$\phi_n=\phi_{P^n}$ defines the required isomorphism, where $H$ is defined as follows
\begin{eqnarray*}
&&\HOM_{\cZ\text{-}\MOD\text{-}\cZ}(\mV_TFP,\mV_TGP)\\
&=&
\HOM_{\cZ\text{-}\MOD\text{-}\cZ}\Big(\HOM_{\mg_T}\big(P,\varprojlim_m FP^m\big),\HOM_{\mg_T}\big(P,\varprojlim_n GP^n\big)\Big)\\
&=&
\varprojlim_n
\HOM_{\cZ\text{-}\MOD\text{-}\cZ}\Big(\varprojlim_m\HOM_{\mg_T}\big(P,FP\big),\HOM_{\mg_T}\big(P,GP\big)\Big)=:H
\end{eqnarray*}
{\it Well-defined:}
Since $\phi$ is a natural transformation, one easily deduces that
  $a^\phi$ is a $\cZ_{T,\Lambda}$-bimodule morphism. For $m\geq n$ we have
  $G(p_{m,n})\circ (a^\phi)_m(\{g_j\}_{j\in\mN})=G(p_{m,n})\circ\phi_m\circ
  g_m=\phi_n\circ F(p_{m,n})\circ g_m=\phi_n\circ
  g_n=(a^\phi)_n(\{g_j\}_{j\in\mN})$ and our map is therefore well-defined.

{\it Injectivity:} Assume $a^\phi=0$, in particular $\phi_n=0$ for any $n\in\mN$. For any tilting
  module $X\in\OTL$ we choose compatible surjections
  $p_n:P^n_{S(0)}(\la)\rightarrow\Delta_{S(0)}(\la)$ for $n\geq 0$. 
Then $F(p_n)$ is again surjective (\cite[Proposition 4.1]{CT}). Since $\phi$ is a natural transformation
  we get $\phi_X\circ F(p_n)=G(p_n)\circ\phi_{n}=0$ for any $n>0$, so
  $\phi_X=0$. Since $\cM_{T,\Lambda}$ is generated by
tilting modules (see Section~\ref{general}) the injectivity follows.   

{\it Surjectivity:} Since (under the assumptions on $T$) any translation through the wall is self-adjoint
  (\cite[Corollary 5.10]{CT}), using Lemma~\ref{freeness} we are allowed to
  assume $F=\ID$.
Let $f\in\HOM_{\cZ\text{-}\MOD\text{-}\cZ}(\mV_TP,\mV_T GP)$. Let
  $X\in\cM_{T,\Lambda}$. We choose a complex
  \begin{eqnarray*}
    \cK: \oplus_{i\in I}P^\infty\rightarrow\oplus_{j\in J}P^\infty\rightarrow
    X\rightarrow 0
  \end{eqnarray*}
with finite sets $I$ and $J$ such that $\mV\cK$ is exact; in other words the homology of $\cK$ does not
contain $L_T(\la)$ as a composition factor. We get a commutative diagram of the
form
\begin{eqnarray*}
  \xymatrix{
  \oplus_{i\in I}\ar@{->}[r]\ar@{->}[d]^{\oplus f}\mV P
& \oplus_{j\in J}\ar@{->}[r]\ar@{->}[d]^{\oplus f}\mV P
& \ar@{->}[r]\ar@{.>}[d]^{f_X}\mV X
& 0\\
 \oplus_{i\in I}\ar@{->}[r]\mV G P
& \oplus_{j\in J}\ar@{->}[r]\mV G P
& \ar@{->}[r]\mV G X
& 0
}
\end{eqnarray*}
where the rows are complexes and where, by assumption,  the first row is
exact. Since $f$ is a $\cZ_{T,\Lambda}$-bimodule map, the left part of the diagram
commutes inducing a unique morphism $f_X$ as indicated. If $X$ is tilting,
then the Structure Theorem~\ref{injectivity} induces a unique map
$f_X\in\HOM_{\mg_T}(X, GX)$ which is natural. By standard arguments, this
induces a natural transformation $\ID\rightarrow G$ when restricted to the
additive category of tilting modules. Since $\cM_{T,\Lambda}$ is generated by
tilting modules the surjectivity follows.
\end{proof}

\subsection{Towards the classification theorem}
We still fix an antidominant block
$\cO_{S_{(0)},\Lambda}$ and denote $R=S(0)$, $\cZ=\cZ_{R,\Lambda}$, and will use the notation from Section~\ref{x}. 
\begin{prop}
\label{l18}
Let $[x]=s_1s_2\cdots s_r$ and $[y]=s_{r+1}s_{r+2}\cdots
s_{r+q}$ be fixed compositions of simple reflections in $\cW(\Lambda)$. Then the following hold 
\begin{enumerate}[a.)] 
  \item \label{aa}
$\HOM_{\MOD\text{-}\cZ}(\cZ_{[x]}(R), \cZ_{[y]}(R))=\HOM_{R\text{-}\MOD{\text-}\cZ}\HOM_Z(Z_{[x]}(R), Z_{[y]}(R))$,
where $R$ acts in the latter by multiplication from the left hand side.
\item \label{bb} 
The space $\HOM_{R\text{-}\MOD\text{-}\cZ}(Z_{[x]}(R), Z_{[y]}(R))$ is a (graded) free left $R$-module of finite rank.
\item The canonical map defines an isomorphism of vector spaces
\small
$$\mC\otimes_R\HOM_{R\text{-}\MOD\text{-}\cZ}(\cZ_{[x]}(R),\cZ_{[y]}(R))
\cong
\HOM_{\mC\otimes_R R\text{-}\MOD\text{-}\cZ}(\cZ_{[x]}(R\otimes_R\mC),
\cZ_{[y]}(R\otimes_R\mC)).$$ 
\normalsize
\item \label{dd}
The canonical map defines an isomorphism 
  \begin{eqnarray*}
&&\cZ\otimes_R\HOM_{R\text{-}\MOD\text{-}\cZ}(\cZ_{[x]}(R), \cZ_{[y]}(R))\\
&\cong&\HOM_{\cZ\text{-}\MOD\text{-}\cZ}(\cZ_{[x]}(\cZ\otimes_RR),
\cZ_{[y]}(\cZ\otimes_RR))\\&=&
\HOM_{\cZ\text{-}\MOD\text{-}\cZ}(\cZ_{[x]}(\cZ), \cZ_{[y]}(\cZ)).
  \end{eqnarray*}
\end{enumerate}
\end{prop}

\begin{proof}
  \begin{enumerate}[a.)]
    \item \label{ZR}
The ring $R$ is obviously
    included in the centre $Z=\cZ_{\cO_{R,\Lambda}}$ via the diagonal embedding
    $R\inj\prod_{\mu\in\cW(\Lambda)} R$, $r\mapsto \{r\}_{\mu}$. The image is
    $s$-invariant for any simple reflection $s$. The claim follows then from
    the definitions. 
\item
We consider the following compositions of translations through walls $F_R=\theta_{s_r}\cdots\theta_{s_2}\theta_{s_1}$,
$G_R=\theta_{s_{r+q}}\cdots\theta_{s_{r+2}}\theta_{s_{r+1}}:\cO_{R,\Lambda}\rightarrow\cO_{R,\Lambda}$
with its specialisations $F_\mC$ and $G_\mC$.
The Structure Theorem~\ref{injectivity} together with \cite[Theorem
4.1]{Fiebigcomb} provides a natural isomorphism of vector spaces, compatible
with the left $R$-action, as follows
\begin{eqnarray*}
\HOM_{\MOD\text{-}\cZ}(\cZ_{[x]}(R), \cZ_{[y]}(R))
\cong\HOM_{\mg_R}(F_R\Delta_R(\la), G_R\Delta_R(\la)).
\end{eqnarray*}
The latter is a free $R$-module of finite rank, because via the tilting
equivalence from Section~\ref{general} we are in the situation of
\cite[Proposition 2.4]{CT} where it is shown that the morphism spaces between
truncated projective modules are free $R$-modules of finite rank. Using part \eqref{ZR} we are done.
\item
We claim that there are natural isomorphisms
\begin{eqnarray*}
&&\mC\otimes_R\HOM_{\MOD\text{-}\cZ}(\cZ_{[x]}(R), \cZ_{[y]}(R))\\
&\stackrel{(1)}\cong&
\mC\otimes_R\HOM_{\mg_R}\big(F_R\Delta_R(\la), G_R\Delta_R(\la)\big)\\
&\stackrel{(2)}\cong&
\HOM_{\mg_\mC}\big(F_\mC\Delta_\mC(\la),G_\mC\Delta_\mC(\la))\\
&\stackrel{(3)}\cong&
\HOM_{\MOD{\text{-}}\mC\otimes_R\cZ}\big(\mV_\mC(F_\mC\Delta_\mC(\la)),\mV_\mC(G_\mC\Delta_\mC(\la)))\\
&\stackrel{(4)}\cong&
\HOM_{\MOD{\text{-}}\mC\otimes_R\cZ}\big(\mV_R(F_R\Delta_R(\la))\otimes_R\mC,\mV_R(G_R\Delta_R(\la))\otimes_R\mC)\\
&\cong&
\HOM_{\mC\otimes_R R\text{-}\MOD\text{-}\cZ}(\cZ_{[x]}(R\otimes_R\mC),
\cZ_{[y]}(R\otimes_R\mC)).
\end{eqnarray*}

The isomorphism $(1)$ follows from the Structure Theorem~\ref{injectivity},
Lemma~\ref{freeness} and Part~\eqref{bb}. The isomorphism $(2)$ exists by \cite[Proposition 2.4, Corollary 5.11]{CT} and the
tilting equivalence. Thanks again to the Structure Theorem and
Corollary~\ref{corcenter} we have
$(3)$. Finally $(4)$ is given by \cite[Proposition 3.11]{Fiebigcomb}.
Using part \eqref{ZR}, the proposition follows.
\item This is completely analogous to \cite[Lemma 12]{SHC}. 
\end{enumerate}
\end{proof}

\subsection{The classification of projective functors}
In this section we formulate and prove the classification theorems for
projective functors. The crucial result is given by the following description of
morphisms between compositions of translation functors as morphisms between
translated deformed Verma modules, where the deformation ring is the centre of
the original category. 
\begin{theorem}
\label{crucial}
Let $T=S_{(0)}$ or $T=\mC$ and let $\cO_{T,\Lambda}$ be an antidominant block
with centre $\cZ_T$. Let $F_T, G_T:\cO_{T,\Lambda}\rightarrow\cO_{T,\Lambda}$ be compositions of
translation through walls with corresponding compositions $F_{\cZ_T},
G_{\cZ_T}:\cO_{\cZ_T,\Lambda}\rightarrow\cO_{\cZ_T,\Lambda}$. 
Then there is an isomorphism of vector spaces (or even of rings if $F=G$)
\begin{eqnarray*}
\HOM(F_T,G_T)\cong\HOM_{\mg_{\cZ_T}}(F_{\cZ_T}\Delta_{\cZ_T}(\la),G_{\cZ_T}\Delta_{\cZ_T}(\la)).
\end{eqnarray*}
\end{theorem}
\begin{proof}
Let first be $T=S_{(0)}=R$. Let $\hat{F}$, $\hat{G}$ denote the functors
  obtained from Lemma~\ref{freeness} such that $\mV_RF\cong\hat{F}\mV_R$ and
  $\mV_RG\cong\hat{G}\mV_R$. Then we have isomorphisms of vector spaces (or rings)
\begin{eqnarray*}
  \begin{array}[t]{ccll}
  &&\HOM(F_R,G_R)\\
&\cong&\HOM_{\cZ\text{-}\MOD\text{-}\cZ}\big(\mV_R F_RP_R^\infty(\la),\mV_R
G_RP_R^\infty(\la)\big)&
\text{(Theorem~\ref{homprojfuncKM})}\\
&\cong&\HOM_{\cZ\text{-}\MOD\text{-}\cZ}\big(\hat{F}_R\mV_R P_R^\infty(\la),\hat{G}_R\mV_R P_R^\infty(\la)\big)\\
&\cong&\HOM_{\cZ\text{-}\MOD\text{-}\cZ}\big(\hat{F}_R (\cZ),\hat{G}_R (\cZ)\big)&
\text{(Theorem~\ref{centerKac})}\\
&\cong&\cZ\otimes_R\HOM_{R\text{-}\MOD\text{-}\cZ}\big(\hat{F}_R (R),\hat{G}_R (R)\big)
&\text{(Proposition~\ref{l18}~\eqref{dd})}
\\
&\cong&\cZ\otimes_R\HOM_{\MOD\text{-}\cZ}\big(\hat{F}_R (R),\hat{G}_R (R)\big)
&\text{(Proposition~\ref{l18}~\eqref{aa})}\\
&\cong&\cZ\otimes_R\HOM_{\MOD\text{-}\cZ}\big(\hat{F}_R\mV_R\Delta_R(\la),\hat{G}_R\mV_R\Delta_R(\la)\big)\\
&\cong&\cZ\otimes_R\HOM_{\MOD\text{-}\cZ}\big(\mV_RF_R\Delta_R(\la),\mV_RG_R\Delta_R(\la)\big)\\
&\cong&
\cZ\otimes_R\HOM_{\mg_R}
\big(F_R\Delta_R(\la), G_R\Delta_R(\la)\big)&\text{(Structure Theorem)}
\\
&\cong&\HOM_{\mg_\cZ}\big(F_R\Delta_\cZ(\la), G_R\Delta_\cZ(\la)\big),
\end{array}
\end{eqnarray*}
the latter by \cite[Proposition 2.4]{CT} and the tilting
  equivalence. Hence the Theorem follows for $T=S_{(0)}$. The case $T=\mC$
  follows then by specialisation or by copying the arguments.  
\end{proof}

We call an endofunctor $F$ of $\OTL$ {\it projective} if it is a direct sum of
direct summands of some compositions of translation through walls.
The following classification of projective functors
(justifying their name) follows immediately:
\begin{cor}[The Classification Theorem]\hfill\\
\label{Class}
Let $T=S_{(0)}$ or $T=\mC$ and let $\cO_{T,\Lambda}$ be an antidominant block
with centre $\cZ_T$. There are natural bijections of isomorphism classes
\begin{eqnarray*}
  \begin{array}[t]{c}
    \big\{\text{indecomposable projective functors on } \cO_{T,\Lambda}\}\\
\updownarrow \text{\small 1:1}\\
\big\{\text{indecomposable tilting objects of }\cO_{T,\Lambda}\big\}\\
\updownarrow \text{\small 1:1}\\
\big\{\text{indecomposable projective objects of }\cO_{T,\tau(\Lambda)}\big\}
  \end{array}
\end{eqnarray*}
\end{cor}

\begin{proof}
  The $F_T\Delta_{\cZ_R}(\la)$'s for $F$ a projective functor are tilting modules in
  $\cO_{\cZ_T,\Lambda}$ (via the tilting equivalence using \cite[Corollary 5.11,
  Proposition 2.4]{CT}). Any tilting module is obtained in such a way
  (\cite[Section 4]{Fiebigcomb}). Moreover, the isomorphism classes of (indecomposable) tilting
  objects in $\cO_{\cZ_T,\Lambda}$ correspond exactly to the ones in
  $\cO_{S_{(0)},\Lambda}$ via specialisation, because $\cZ$ is local. Hence the classification follows. 
\end{proof}

We have a Krull-Remak-Schmidt property for projective functors:

\begin{cor}
\label{decomp}
  Let $F:\cO_{\mC,\Lambda}\rightarrow\cO_{\mC,\Lambda}$ be a projective
  functor. Then
\begin{enumerate}[a.)]
  \item $F\cong\bigoplus_{i=1}^r F_i$ for some indecomposable functors $F_i$ and $r\in
  \mN$.
\item Moreover, $F(\Delta_\mC(\lambda))\cong\bigoplus_{i=1}^r F_i(\Delta_\mC(\la))$ is a decomposition into indecomposable direct summands.
\end{enumerate}
\end{cor}
\begin{proof}
  The claim follows directly from the previous Theorem~\ref{crucial} and Corollary~\ref{Class}.
\end{proof}

\bibliography{ref2}

\noindent C. Stroppel, Department of Mathematics, University of Glasgow,
15 University Gardens, Glasgow G12 8QW, (United Kingdom).
\end{document}